\documentclass[reqno,10pt,centertags,draft]{amsart}
\usepackage{amsmath,amsthm,amscd,amssymb,latexsym,esint,upref,stmaryrd,
enumerate,color,verbatim,yfonts}
\usepackage{graphicx}
\usepackage{hyperref}
\newcommand*{\mailto}[1]{\href{mailto:#1}{\nolinkurl{#1}}}
\newcommand{\arxiv}[1]{\href{http://arxiv.org/abs/#1}{arXiv:#1}}

\makeatletter
\def\widebreve{\mathpalette\wide@breve}
\def\wide@breve#1#2{\sbox\z@{$#1#2$}%
     \mathop{\vbox{\m@th\ialign{##\crcr
\kern0.08em\brevefill#1{0.8\wd\z@}\crcr\noalign{\nointerlineskip}%
                    $\hss#1#2\hss$\crcr}}}\limits}
\def\brevefill#1#2{$\m@th\sbox\tw@{$#1($}%
  \hss\resizebox{#2}{\wd\tw@}{\rotatebox[origin=c]{90}{\upshape(}}\hss$}
\makeatletter

\newcommand{\R}{{\mathbb R}}
\newcommand{\N}{{\mathbb N}}
\newcommand{\Z}{{\mathbb Z}}

\newcommand{\bbC}{{\mathbb{C}}}

\newcommand{\bbN}{{\mathbb{N}}}

\newcommand{\bbR}{{\mathbb{R}}}

\newcommand{\bbZ}{{\mathbb{Z}}}

\newcommand{\cB}{{\mathcal B}}

\newcommand{\cH}{{\mathcal H}}

\newcommand{\cN}{{\mathcal N}}

\newcommand{\beq}{\begin{equation}}
\newcommand{\enq}{\end{equation}}

\renewcommand{\a}{\alpha}
\renewcommand{\b}{\beta}
\newcommand{\g}{\gamma}
\renewcommand{\d}{\delta}

\renewcommand{\l}{\lambda}

\newcommand{\s}{\sigma}

\newcommand{\G}{\Gamma}

\DeclareMathOperator{\supp}{supp}

\DeclareMathOperator{\dom}{dom}

\renewcommand{\Re}{\text{\rm Re}}
\renewcommand{\Im}{\text{\rm Im}}
\renewcommand{\ln}{\text{\rm ln}}

\newcommand{\no}{\notag}
\newcommand{\lb}{\label}
\newcommand{\f}{\frac}

\newcommand{\ol}{\overline}
\newcommand{\bs}{\backslash}

\newcommand{\wti}{\widetilde}
\newcommand{\Oh}{O}
\newcommand{\oh}{o}
\newcommand{\hatt}{\widehat} 
\newcommand{\dott}{\,\cdot\,}
\renewcommand{\Dot}{\overset{\textbf{\Large.}}}

\newcommand{\dotA}{{\hspace{0.08cm}\Dot{\hspace{-0.08cm} A}}}

\newcommand{\bi}{\bibitem}

\let\geq\geqslant
\let\leq\leqslant

\newcommand{\lam}{\lambda}

\newcommand{\al}{\alpha}

\newcommand{\ACl}{{AC_{loc}((a,b))}}
\newcommand{\Ll}{{L^1_{loc}((a,b);dx)}}
\newcommand{\SL}{\text{\textnormal{SL}}}

\newcommand{\high}[1]{{\raisebox{1mm}{$#1$}}}

\makeatletter
\def\theequation{\@arabic\c@equation}

\allowdisplaybreaks 
\numberwithin{equation}{section}

\newtheorem{theorem}{Theorem}[section]

\newtheorem{definition}[theorem]{Definition}
\newtheorem{hypothesis}[theorem]{Hypothesis}

\theoremstyle{remark}
\newtheorem{remark}[theorem]{Remark}

\begin{document}

\title[Jacobi Donoghue $m$-functions]{The Jacobi Operator and its Donoghue $m$-Functions}

\author[F.\ Gesztesy]{Fritz Gesztesy}
\address{Department of Mathematics, 
Baylor University, Sid Richardson Bldg., 1410 S.\,4th Street, Waco, TX 76706, USA}
\email{\mailto{Fritz\_Gesztesy@baylor.edu}}
\urladdr{\url{http://www.baylor.edu/math/index.php?id=935340}}

\author[M. Piorkowski]{Mateusz Piorkowski}
\address{Mathematical Sciences and Research Institute, UC Berkeley, 17 Gauss Way, CA 94720, USA}
\email{\href{mailto:mathpiorkowski@gmail.com}{Mathpiorkowski@gmail.com}}

\author[J.\ Stanfill]{Jonathan Stanfill}
\address{Department of Mathematics, 
Baylor University, Sid Richardson Bldg., 1410 S.\,4th Street, Waco, TX 76706, USA}
\email{\mailto{Jonathan\_Stanfill@baylor.edu}}
\urladdr{\url{http://sites.baylor.edu/jonathan-stanfill/}}


\date{\today}
\@namedef{subjclassname@2020}{\textup{2020} Mathematics Subject Classification}

\begin{abstract} 
In this paper we construct Donoghue $m$-functions for the Jacobi differential operator in $L^2\big((-1,1); (1-x)^\a (1+x)^\b dx\big)$, associated to the differential expression
\begin{align*}
\begin{split} 
\tau_{\a,\b} = - (1-x)^{-\a} (1+x)^{-\b}(d/dx) \big((1-x)^{\a+1}(1+x)^{\b+1}\big) (d/dx),&     \\ 
x \in (-1,1), \; \a, \b \in \bbR,
\end{split} 
\end{align*}
whenever at least one endpoint, $x=\pm1$, is in the limit circle case.
In doing so, we provide a full treatment of the Jacobi operator's $m$-functions corresponding to coupled boundary conditions whenever both endpoints are in the limit circle case, a topic not covered in the literature.
\end{abstract}

\maketitle

{\scriptsize{\tableofcontents}}

\maketitle



\section{Introduction} \lb{s1}

This paper should be regarded as a sequel to the recent \cite{GLNPS21} in which the Donoghue $m$-function was derived for singular Sturm--Liouville operators. To illustrate the theory, we now apply it to a representative example, the Jacobi differential operator associated with $L^2\big((-1,1); (1-x)^\a (1+x)^\b dx\big)$-realizations of the the differential expression,
\begin{align} \lb{1.0}
\begin{split} 
\tau_{\a,\b} = - (1-x)^{-\a} (1+x)^{-\b}(d/dx) \big((1-x)^{\a+1}(1+x)^{\b+1}\big) (d/dx),&     \\ 
x \in (-1,1), \; \a, \b \in \bbR, 
\end{split} 
\end{align}
whenever at least one endpoint, $x=\pm1$, is in the limit circle case (see, e.g. \cite[Ch.~22]{AS72}, \cite{BEZ01}, \cite{Bo29}, \cite[Sect.~23]{Ev05}, \cite[Ch.~4]{Is05}, \cite[Sects.~VII.6.1, XIV.2]{Kr02}, \cite[Ch.~18]{Ov20}, \cite[Ch.~IV]{Sz75}). In particular, this provides a full treatment of $m$-functions corresponding to coupled boundary conditions whenever both endpoints are in the limit circle case, a new result.

To set the stage we briefly discuss abstract Donoghue $m$-functions following \cite{GNWZ19}, \cite{GKMT01}, \cite{GLNPS21}, and \cite{GMT98}. Given a self-adjoint extension $A$ of a densely defined, closed, symmetric operator $\dotA$ in $\cH$ (a complex, separable Hilbert space) with equal deficiency indices and the deficiency subspace $\cN_i$ of $\dotA$ in $\cH$, with
\begin{equation}
 \cN_i = \ker \big(\big(\dotA\big)\high{^*} - i I_{\cH}\big), \quad
\dim \, (\cN_i)=k\in \bbN \cup \{\infty\},     \lb{1.1}
\end{equation}
the Donoghue $m$-operator $M_{A,\cN_i}^{Do} (\dott) \in\cB(\cN_i)$ associated with the pair
 $(A,\cN_i)$  is given by
\begin{align}
\begin{split}
M_{A,\cN_i}^{Do}(z)&=P_{\cN_i} (zA + I_\cH)(A - z I_{\cH})^{-1}
P_{\cN_i}\big\vert_{\cN_i}      \\
&=zI_{\cN_i} + (z^2+1) P_{\cN_i} (A - z I_{\cH})^{-1}
P_{\cN_i} \big\vert_{\cN_i}\,, \quad  z\in \bbC\backslash \bbR,     \lb{1.2}
\end{split}
\end{align}
with $I_{\cN_i}$ the identity operator in $\cN_i$, and $P_{\cN_i}$ the orthogonal projection in
$\cH$ onto $\cN_i$. The special case $k=1$, was discussed in detail by Donoghue \cite{Do65}; for the case $k\in \bbN$ we refer to \cite{GT00}.

More generally, given a self-adjoint extension $A$ of $\dotA$ in $\cH$ and a closed,
linear subspace $\cN$ of $\cN_i$,
the Donoghue $m$-operator $M_{A,\cN}^{Do} (\dott) \in\cB(\cN)$ associated with the pair
 $(A,\cN)$  is defined by
\begin{align}
\begin{split}
M_{A,\cN}^{Do}(z)&=P_\cN (zA + I_\cH)(A - z I_{\cH})^{-1}
P_\cN\big\vert_\cN      \\
&=zI_\cN+(z^2+1)P_\cN(A - z I_{\cH})^{-1}
P_\cN\big\vert_\cN\,, \quad  z\in \bbC\backslash \bbR,     \lb{1.3}
\end{split}
\end{align}
with $I_\cN$ the identity operator in $\cN$ and $P_\cN$ the orthogonal projection in $\cH$
onto $\cN$. 

Since $M_{A,\cN}^{Do}(z)$ is analytic for $z\in \bbC\backslash\bbR$ and satisfies (see \cite[Theorem~5.3]{GNWZ19}) 
\begin{align}
& [\Im(z)]^{-1} \Im\big(M_{A,\cN}^{Do} (z)\big) \geq
2 \Big[\big(|z|^2 + 1\big) + \big[\big(|z|^2 -1\big)^2 + 4 (\Re(z))^2\big]^{1/2}\Big]^{-1} I_{\cN},    \no \\
& \hspace*{9.5cm}  z\in \bbC\backslash \bbR,    \lb{1.4}
\end{align}
$M_{A,\cN}^{Do}(\dott)$ is a $\cB(\cN)$-valued Nevanlinna--Herglotz function. Thus, $M_{A,\cN}^{Do}(\dott)$ 
admits the representation
\begin{equation}
M_{A,\cN}^{Do}(z) = \int_\bbR
d\Omega_{A,\cN}^{Do}(\lambda) \bigg[\f{1}{\lambda-z} -
\f{\lambda}{\lambda^2 + 1}\bigg], \quad z\in\bbC\backslash\bbR,    \lb{1.5}
\end{equation}
where the $\cB(\cN)$-valued measure $\Omega_{A,\cN}^{Do}(\dott)$ satisfies
\begin{align}
&\Omega_{A,\cN}^{Do}(\lambda)=(\lambda^2 + 1) (P_{\cN} E_A(\lambda)P_{\cN}\big\vert_{\cN}),
\lb{1.6} \\
&\int_\bbR d\Omega_{A,\cN}^{Do}(\lambda) \, (1+\lambda^2)^{-1}=I_{\cN},    \lb{1.7} \\
&\int_\bbR d(\xi,\Omega_{A,\cN}^{Do} (\lambda)\xi)_{\cN} = \infty \, \text{ for all } \,
 \xi \in \cN \backslash \{0\},    \lb{1.8}
\end{align}
with $E_A(\dott)$ the family of strongly right-continuous spectral projections of $A$ in $\cH$ (see \cite{GKMT01} for details). Operators of the type $M_{A,\cN}^{Do}(\dott)$ and some of its variants have attracted
considerable attention in the literature. They appear to go back to Krein \cite{Kr46} (see also 
\cite{KL71}), Saakjan \cite{Sa65}, and independently, Donoghue \cite{Do65}. The interested reader can find a wealth of additional information in the context of \eqref{1.2}--\eqref{1.8} in 
 \cite{AB09}, \cite{AP04}, \cite{BL07}--\cite{BR16}, \cite{BMN02}, \cite{BGW09}--\cite{BGP08},
\cite{DHMdS09}--\cite{DMT88}, \cite{GKMT01}--\cite{GT00}, \cite{GWZ13b}, \cite{HMM13}, \cite{KO78}, \cite{LT77}, \cite{Ma92a}, \cite{MN11}, \cite{MN12}, \cite{Ma04}, \cite{Mo09}, \cite{Na89}--\cite{Na94}, \cite{Pa13}, \cite{Po04}, \cite{Ry07}, and the references therein.

Without going into further details (see \cite[Corollary~5.8]{GNWZ19} for details) we note that the prime reason for the interest in $M_{A,\cN_i}^{Do}(\dott)$ lies in the fundamental fact that the entire spectral information of $A$ contained in its family of spectral projections $E_A(\dott)$, is already encoded in the $\cB(\cN_i)$-valued measure $\Omega_{A,\cN_i}^{Do}(\dott)$  (including multiplicity properties of the spectrum of $A$) if and only if $\dotA$ is completely non-self-adjoint in $\cH$ (that is, if and only if $\dotA$ has no invariant subspace on which it is self-adjoint, see 
\cite[Lemma~5.4]{GNWZ19}). 

We also note that a particularly attractive feature of the Donoghue $m$-operator, that distinguishes it from the 
Weyl--Titchmarsh--Kodaira $m$-operator, consists of the explicit appearance of the resolvent $(A - z I_{\cH})^{-1}$, 
$z \in \bbC \backslash \bbR$, in it's definition \eqref{1.2} (resp., \eqref{1.3}).

In the remainder of this paper, we will exclusively focus on the particular case
$\cN = \cN_i = \ker\big(\big(\dotA\big)\high{^*} - i I_{\cH}\big)$, with $\dotA$ being a singular Sturm--Liouville operator.

Turning to the content of each section, we discuss the necessary background in connection to singular Sturm--Liouville operators in Section \ref{s2}. In Sections \ref{s4} and \ref{s3} we recall the Donoghue $m$-functions in the two limit circle and one limit circle endpoint cases, respectively, following \cite{GLNPS21}. The Jacobi operator and its Donoghue $m$-functions are the topic of Section \ref{s6}, with Appendices \ref{sA}--\ref{sC} providing a detailed treatment of solutions of the Jacobi differential equation and the associated hypergeometric differential equations.

Finally, some comments on some of the basic notation used throughout this paper.  If $T$ is a linear operator mapping (a subspace of) a Hilbert space into another, then $\dom(T)$ and $\ker(T)$ denote the domain and kernel (i.e., null space) of $T$. The spectrum and resolvent set of a closed linear operator in a Hilbert space will be denoted by $\sigma(\dott)$ and $\rho(\dott)$, respectively. Moreover, we denote the scalar product and norm in  $L^2((a,b);rdx)$ by 
$(\dott,\dott)_{L^2((a,b);rdx)}$ (linear in the second argument) and $\|\dott\|_{L^2((a,b);rdx)}$.

\section{Some Background} \lb{s2}

In this section we briefly recall the basics of singular Sturm--Liouville operators. The material is standard and can be found, for instance, in \cite[Ch.~6]{BHS20}, \cite[Chs.~8, 9]{CL85}, \cite[Sects.~13.6, 13.9, 13.10]{DS88}, \cite{EGNT13}, \cite[Ch.~4]{GZ21}, \cite[Ch.~III]{JR76}, \cite[Ch.~V]{Na68}, \cite{NZ92}, \cite[Ch.~6]{Pe88}, \cite[Ch.~9]{Te14}, \cite[Sect.~8.3]{We80}, \cite[Ch.~13]{We03}, \cite[Chs.~4, 6--8]{Ze05}.

Throughout this section we make the following assumptions:

\begin{hypothesis} \lb{h2.1}
Let $(a,b) \subseteq \bbR$ and suppose that $p,q,r$ are $($Lebesgue\,$)$ measurable functions on $(a,b)$ 
such that the following items $(i)$--$(iii)$ hold: \\[1mm] 
$(i)$ \hspace*{1.1mm} $r>0$ a.e.~on $(a,b)$, $r\in\Ll$. \\[1mm] 
$(ii)$ \hspace*{.1mm} $p>0$ a.e.~on $(a,b)$, $1/p \in\Ll$. \\[1mm] 
$(iii)$ $q$ is real-valued a.e.~on $(a,b)$, $q\in\Ll$. 
\end{hypothesis}

Given Hypothesis \ref{h2.1}, we study Sturm--Liouville operators associated with the general, 
three-coefficient differential expression $\tau$ of the form
\begin{equation}
\tau=\f{1}{r(x)}\left[-\f{d}{dx}p(x)\f{d}{dx} + q(x)\right] \, \text{ for a.e.~$x\in(a,b) \subseteq \bbR$.}     \lb{2.1}
\end{equation} 

If $f\in AC_{loc}((a,b))$, then the quasi-derivative of $f$ is defined to be $f^{[1]}:=pf'$.  Moreover, the Wronskian of two functions $f,g\in AC_{loc}((a,b))$ is defined by
\begin{equation}
W(f,g)(x) = f(x)g^{[1]}(x)-f^{[1]}(x)g(x)\,\text{ for a.e.~$x\in (a,b)$}.
\end{equation}

Assuming Hypothesis \ref{h2.1}, the {\it maximal operator} $T_{max}$ in $L^2((a,b);rdx)$ associated with $\tau$ is defined by
\begin{align}
&T_{max} f = \tau f,    \no
\\
& f \in \dom(T_{max})=\big\{g \in L^2((a,b);rdx) \, \big| \,g,g^{[1]}\in\ACl;   \lb{2.2} \\ 
& \hspace*{6.3cm}  \tau g \in L^2((a,b);rdx)\big\}.   \no
\end{align}
The {\it preminimal operator} $\dot T_{min} $ in $L^2((a,b);rdx)$ associated with $\tau$ is defined by 
\begin{align}
&\dot T_{min}  f = \tau f,   \no
\\
&f \in \dom \big(\dot T_{min}\big)=\big\{g \in L^2((a,b);rdx) \, \big| \, g,g^{[1]}\in\ACl;   \lb{2.3}
\\
&\hspace*{3.1cm} \supp \, (g)\subset(a,b) \text{ is compact; } \tau g \in L^2((a,b);rdx)\big\}.   \no
\end{align}
One can prove that $\dot T_{min} $ is closable, and one then defines the {\it minimal operator} $T_{min}$ by 
$T_{min} = \ol{\dot T_{min}}$.

Still assuming Hypothesis \ref{h2.1}, one can prove the following basic fact, 
\begin{equation} 
\big(\dot T_{min}\big)^* = T_{max},  
\end{equation} 
and hence $T_{max}$ is closed. Moreover, $\dot T_{min} $ is essentially self-adjoint if and only if $T_{max}$ is symmetric, and then $\ol{\dot T_{min} }=T_{min}=T_{max}$.

The celebrated Weyl alternative can be stated as follows:

\begin{theorem}[Weyl's Alternative] \lb{t2.5} ${}$ \\
Assume Hypothesis \ref{h2.1}. Then the following alternative holds: Either \\[1mm] 
$(i)$ for every $z\in\bbC$, all solutions $u$ of $(\tau-z)u=0$ are in $L^2((a,b);rdx)$ near $b$ 
$($resp., near $a$$)$, \\[1mm] 
or, \\[1mm] 
$(ii)$ for every $z\in\bbC$, there exists at least one solution $u$ of $(\tau-z)u=0$ which is not in $L^2((a,b);rdx)$ near $b$ $($resp., near $a$$)$. In this case, for each $z\in\bbC\bs\bbR$, there exists precisely one solution $u_b$ $($resp., $u_a$$)$ of $(\tau-z)u=0$ $($up to constant multiples$)$ which lies in $L^2((a,b);rdx)$ near $b$ $($resp., near $a$$)$. 
\end{theorem}

This yields the limit circle/limit point classification of $\tau$ at an interval endpoint as follows. 

\begin{definition} \lb{d2.6} 
Assume Hypothesis \ref{h2.1}. \\[1mm]  
In case $(i)$ in Theorem \ref{t2.5}, $\tau$ is said to be in the {\it limit circle case} at $b$ $($resp., $a$$)$. $($Frequently, $\tau$ is then called {\it quasi-regular} at $b$ $($resp., $a$$)$.$)$
\\[1mm] 
In case $(ii)$ in Theorem \ref{t2.5}, $\tau$ is said to be in the {\it limit point case} at $b$ $($resp., $a$$)$. \\[1mm]
If $\tau$ is in the limit circle case at $a$ and $b$ then $\tau$ is also called {\it quasi-regular} on $(a,b)$. 
\end{definition}

The next result links self-adjointness of $T_{min}$ (resp., $T_{max}$) and the limit point property of $\tau$ at both endpoints.  Here, and throughout, we shall employ the notation
\begin{equation}
\cN_z = \ker\big(T_{max}-zI_{L^2((a,b);rdx)}\big),\quad z\in \bbC.
\end{equation}

\begin{theorem} \lb{t2.7}
Assume Hypothesis~\ref{h2.1}, then the following items $(i)$ and $(ii)$ hold: \\[1mm] 
$(i)$ If $\tau$ is in the limit point case at $a$ $($resp., $b$$)$, then 
\begin{equation} 
W(f,g)(a)=0 \, \text{$($resp., $W(f,g)(b)=0$$)$ for all $f,g\in\dom(T_{max})$.} 
\end{equation} 
$(ii)$ Let $T_{min}=\ol{\dot T_{min} }$. Then
\begin{align}
\begin{split}
n_\pm(T_{min}) &= \dim\big(\cN_{\pm i}\big)    \\
& = \begin{cases}
2 & \text{if $\tau$ is in the limit circle case at $a$ and $b$,}\\
1 & \text{if $\tau$ is in the limit circle case at $a$} \\
& \text{and in the limit point case at $b$, or vice versa,}\\
0 & \text{if $\tau$ is in the limit point case at $a$ and $b$}.
\end{cases}
\end{split} 
\end{align}
In particular, $T_{min} = T_{max}$ is self-adjoint if and only if $\tau$ is in the limit point case at $a$ and $b$. 
\end{theorem}

All self-adjoint extensions of $T_{min}$ are then described as follows:

\begin{theorem} \lb{t2.8}
Assume Hypothesis \ref{h2.1} and that $\tau$ is in the limit circle case at $a$ and $b$ $($i.e., $\tau$ is quasi-regular 
on $(a,b)$$)$. In addition, assume that 
$v_j \in \dom(T_{max})$, $j=1,2$, satisfy 
\begin{equation}
W(\ol{v_1}, v_2)(a) = W(\ol{v_1}, v_2)(b) = 1, \quad W(\ol{v_j}, v_j)(a) = W(\ol{v_j}, v_j)(b) = 0, \; j= 1,2.  
\end{equation}
$($E.g., real-valued solutions $v_j$, $j=1,2$, of $(\tau - \lambda) u = 0$ with $\lambda \in \bbR$, such that 
$W(v_1,v_2) = 1$.$)$ For $g\in\dom(T_{max})$ we introduce the generalized boundary values 
\begin{align}
\begin{split} 
\wti g_1(a) &= - W(v_2, g)(a), \quad \wti g_1(b) = - W(v_2, g)(b),    \\
\wti g_2(a) &= W(v_1, g)(a), \quad \;\,\,\, \wti g_2(b) = W(v_1, g)(b).   \lb{2.10}
\end{split} 
\end{align}
Then the following items $(i)$--$(iii)$ hold: \\[1mm] 
$(i)$ All self-adjoint extensions $T_{\g,\d}$ of $T_{min}$ with separated boundary conditions are of the form
\begin{align}
& T_{\g,\d} f = \tau f, \quad \g,\d\in[0,\pi),  \lb{2.11} \\
&f\in \dom(T_{\g,\d}) = \bigg\{g\in \dom(T_{max})\,\bigg|\,
\begin{aligned}
\cos(\g)\wti g_1(a)+\sin(\g)\wti g_2(a)&=0,\no\\
\cos(\d)\wti g_1(b)+\sin(\d)\wti g_2(b)&=0
\end{aligned}
\bigg\}.\no 
\end{align}
$(ii)$ All self-adjoint extensions $T_{\varphi,R}$ of $T_{min}$ with coupled boundary conditions are of the type
\begin{align}
\begin{split} 
& T_{\varphi,R} f = \tau f,    \\
& f \in \dom(T_{\varphi,R})=\bigg\{g\in\dom(T_{max}) \, \bigg| \begin{pmatrix} \wti g_1(b)\\ \wti g_2(b)\end{pmatrix} 
= e^{i\varphi}R \begin{pmatrix}
\wti g_1(a)\\ \wti g_2(a)\end{pmatrix} \bigg\}, \lb{2.12}
\end{split}
\end{align}
where $\varphi\in[0,\pi)$, and $R$ is a real $2\times2$ matrix with $\det(R)=1$ 
$($i.e., $R \in SL(2,\bbR)$$)$.  \\[1mm] 
$(iii)$ Every self-adjoint extension of $T_{min}$ is either of type $(i)$ $($i.e., separated\,$)$ or of type 
$(ii)$ $($i.e., coupled\,$)$.
\end{theorem}

One can now detail the characterization of $\dom(T_{min})$ by
\begin{align}
\begin{split} 
&T_{min} f = \tau f,      \lb{2.17} \\
&f \in \dom(T_{min}) = \big\{g\in\dom(T_{max})  \, \big| \, \wti g_1 (a) = \wti g_2 (a) = \wti g_1 (b) = \wti g_2 (b) = 0\big\}.   \end{split} 
\end{align} 

Finally, we turn to the characterization of generalized boundary values in the case where $T_{min}$ is bounded from below following \cite{GLN20} and \cite{NZ92}.

We briefly recall the basics of oscillation theory with particular emphasis on principal and nonprincipal solutions, a notion originally due to Leighton and Morse \cite{LM36} (see also Rellich \cite{Re43}, \cite{Re51} and Hartman and Wintner \cite[Appendix]{HW55}). Our outline below follows \cite{CGN16}, 
\cite[Sects.~13.6, 13.9, 13.10]{DS88}, \cite[Ch.~7]{GZ21}, \cite[Ch.~XI]{Ha02}, \cite{NZ92}, \cite[Chs.~4, 6--8]{Ze05}. 

\begin{definition} \lb{d2.10}
Assume Hypothesis \ref{h2.1}. \\[1mm] 
$(i)$ Fix $c\in (a,b)$ and $\lambda\in\bbR$. Then $\tau - \lam$ is
called {\it nonoscillatory} at $a$ $($resp., $b$$)$, 
if every real-valued solution $u(\lambda,\dott)$ of 
$\tau u = \lambda u$ has finitely many
zeros in $(a,c)$ $($resp., $(c,b)$$)$. Otherwise, $\tau - \lam$ is called {\it oscillatory}
at $a$ $($resp., $b$$)$. \\[1mm] 
$(ii)$ Let $\lambda_0 \in \bbR$. Then $T_{min}$ is called bounded from below by $\lambda_0$, 
and one writes $T_{min} \geq \lambda_0 I_{L^2((a,b);rdx)}$, if 
\begin{equation} 
\big(u, [T_{min} - \lambda_0 I_{L^2((a,b);rdx)}]u\big)_{L^2((a,b);rdx)}\geq 0, \quad u \in \dom(T_{min}).
\end{equation}
\end{definition}

The following is a key result. 

\begin{theorem} \lb{t2.11} 
Assume Hypothesis \ref{h2.1}. Then the following items $(i)$--$(iii)$ are
equivalent: \\[1mm] 
$(i)$ $T_{min}$ $($and hence any symmetric extension of $T_{min})$
is bounded from below. \\[1mm] 
$(ii)$ There exists a $\nu_0\in\bbR$ such that for all $\lambda < \nu_0$, $\tau - \lam$ is
nonoscillatory at $a$ and $b$. \\[1mm]
$(iii)$ For fixed $c, d \in (a,b)$, $c \leq d$, there exists a $\nu_0\in\bbR$ such that for all
$\lambda<\nu_0$, $\tau u = \lambda u$ has $($real-valued\,$)$ nonvanishing solutions
$u_a(\lambda,\dott) \neq 0$,
$\hatt u_a(\lambda,\dott) \neq 0$ in the neighborhood $(a,c]$ of $a$, and $($real-valued\,$)$ nonvanishing solutions
$u_b(\lambda,\dott) \neq 0$, $\hatt u_b(\lambda,\dott) \neq 0$ in the neighborhood $[d,b)$ of
$b$, such that 
\begin{align}
&W(\hatt u_a (\lambda,\dott),u_a (\lambda,\dott)) = 1,
\quad u_a (\lambda,x)=\oh(\hatt u_a (\lambda,x))
\text{ as $x\downarrow a$,} \lb{2.18} \\
&W(\hatt u_b (\lambda,\dott),u_b (\lambda,\dott))\, = 1,
\quad u_b (\lambda,x)\,=\oh(\hatt u_b (\lambda,x))
\text{ as $x\uparrow b$,} \lb{2.19} \\
&\int_a^c dx \, p(x)^{-1}u_a(\lambda,x)^{-2}=\int_d^b dx \, 
p(x)^{-1}u_b(\lambda,x)^{-2}=\infty,  \lb{2.20} \\
&\int_a^c dx \, p(x)^{-1}{\hatt u_a(\lambda,x)}^{-2}<\infty, \quad 
\int_d^b dx \, p(x)^{-1}{\hatt u_b(\lambda,x)}^{-2}<\infty. \lb{2.21}
\end{align}
\end{theorem}

\begin{definition} \lb{d2.12}
Assume Hypothesis \ref{h2.1}, suppose that $T_{min}$ is bounded from below, and let 
$\lambda\in\bbR$. Then $u_a(\lambda,\dott)$ $($resp., $u_b(\lambda,\dott)$$)$ in Theorem
\ref{t2.11}\,$(iii)$ is called a {\it principal} $($or {\it minimal}\,$)$
solution of $\tau u=\lambda u$ at $a$ $($resp., $b$$)$. A real-valued solution 
${\widebreve u}_a(\lambda,\dott)$ $($resp., ${\widebreve u}_b(\lambda,\dott)$$)$ of $\tau
u=\lambda u$ linearly independent of $u_a(\lambda,\dott)$ $($resp.,
$u_b(\lambda,\dott)$$)$ is called {\it nonprincipal} at $a$ $($resp., $b$$)$. In particular, $\hatt u_a (\lambda,\dott)$ 
$($resp., $\hatt u_b(\lambda,\dott)$$)$ in \eqref{2.18}--\eqref{2.21} are nonprincipal solutions at $a$ $($resp., $b$$)$.  
\end{definition}

Next, we revisit  in Theorem \ref{t2.8} how the generalized boundary values are utilized in the description of all self-adjoint extensions of $T_{min}$ in the case where $T_{min}$ is bounded from below. 

\begin{theorem} [{\cite[Theorem~4.5]{GLN20}}]\lb{t2.13}
Assume Hypothesis \ref{h2.1} and that $\tau$ is in the limit circle case at $a$ and $b$ $($i.e., $\tau$ is quasi-regular 
on $(a,b)$$)$. In addition, assume that $T_{min} \geq \lambda_0 I$ for some $\lambda_0 \in \bbR$, and denote by 
$u_a(\lambda_0, \dott)$ and $\hatt u_a(\lambda_0, \dott)$ $($resp., $u_b(\lambda_0, \dott)$ and 
$\hatt u_b(\lambda_0, \dott)$$)$ principal and nonprincipal solutions of $\tau u = \lambda_0 u$ at $a$ 
$($resp., $b$$)$, satisfying
\begin{equation}
W(\hatt u_a(\lambda_0,\dott), u_a(\lambda_0,\dott)) = W(\hatt u_b(\lambda_0,\dott), u_b(\lambda_0,\dott)) = 1.  
\lb{2.22} 
\end{equation}
Introducing $v_j \in \dom(T_{max})$, $j=1,2$, via 
\begin{align}
v_1(x) = \begin{cases} \hatt u_a(\lambda_0,x), & \text{for $x$ near a}, \\
\hatt u_b(\lambda_0,x), & \text{for $x$ near b},  \end{cases}   \quad 
v_2(x) = \begin{cases} u_a(\lambda_0,x), & \text{for $x$ near a}, \\
u_b(\lambda_0,x), & \text{for $x$ near b},  \end{cases}   \lb{2.23}
\end{align} 
one obtains for all $g \in \dom(T_{max})$, 
\begin{align}
\begin{split} 
\wti g(a) &= - W(v_2, g)(a) = \wti g_1(a) =  - W(u_a(\lambda_0,\dott), g)(a) = \lim_{x \downarrow a} \f{g(x)}{\hatt u_a(\lambda_0,x)},    \\
\wti g(b) &= - W(v_2, g)(b) = \wti g_1(b) =  - W(u_b(\lambda_0,\dott), g)(b) 
= \lim_{x \uparrow b} \f{g(x)}{\hatt u_b(\lambda_0,x)},    \lb{2.24} 
\end{split} \\
\begin{split} 
{\wti g}^{\, \prime}(a) &= W(v_1, g)(a) = \wti g_2(a) = W(\hatt u_a(\lambda_0,\dott), g)(a) 
= \lim_{x \downarrow a} \f{g(x) - \wti g(a) \hatt u_a(\lambda_0,x)}{u_a(\lambda_0,x)},    \\ 
{\wti g}^{\, \prime}(b) &= W(v_1, g)(b) = \wti g_2(b) = W(\hatt u_b(\lambda_0,\dott), g)(b)  
= \lim_{x \uparrow b} \f{g(x) - \wti g(b) \hatt u_b(\lambda_0,x)}{u_b(\lambda_0,x)}.    \lb{2.25}
\end{split} 
\end{align}
In particular, the limits on the right-hand sides in \eqref{2.24}, \eqref{2.25} exist. 
\end{theorem}

The Friedrichs extension $T_F$ of $T_{min}$ now permits a particularly simple characterization in terms of the generalized boundary values $\wti g(a), \wti g(b)$ as derived by Niessen and Zettl \cite{NZ92}(see also \cite{GP79}, \cite{Ka72}, \cite{Ka78}, \cite{KKZ86}, \cite{MZ00}, \cite{Re51}, \cite{Ro85}, \cite{YSZ15}):

\begin{theorem} \lb{t2.15}
Assume Hypothesis \ref{h2.1} and that $\tau$ is in the limit circle case at $a$ and $b$ $($i.e., $\tau$ is quasi-regular 
on $(a,b)$$)$. In addition, assume that $T_{min} \geq \lambda_0 I$ for some $\lambda_0 \in \bbR$. Then the Friedrichs extension $T_F=T_{0,0}$ of $T_{min}$ is characterized by
\begin{align}
T_F f = \tau f, \quad f \in \dom(T_F)= \big\{g\in\dom(T_{max})  \, \big| \, \wti g(a) = \wti g(b) = 0\big\}.    \lb{2.26}
\end{align}
\end{theorem}

\section{Donoghue $m$-functions: Two Limit Circle Endpoints} \lb{s4}

The Donoghue $m$-functions in the case where $\tau$ is in the limit circle case at $a$ and $b$ is the primary topic of this section following \cite[Sect. 6]{GLNPS21}. 

\begin{hypothesis}\lb{h4.1}
In addition to Hypothesis \ref{h2.1} assume that $\tau$ is in the limit circle case at $a$ and $b$. Moreover, for $z\in \rho(T_{0,0})$, let $\{u_j(z,\dott)\}_{j=1,2}$ denote solutions to $\tau u = zu$ which satisfy the boundary conditions
\begin{equation}\lb{3.1}
\begin{split}
\wti u_1(z,a)=0, &\quad \wti u_1(z,b)=1,\\
\wti u_2(z,a)=1, &\quad \wti u_2(z,b)=0.
\end{split}
\end{equation}
\end{hypothesis}

Assume Hypotheses \ref{h4.1}. By Theorem \ref{t2.8} or Theorem \ref{t2.13}, the following statements $(i)$--$(iii)$ hold.\\[1mm]
$(i)$  If $\g,\d\in [0,\pi)$, then the operator $T_{\g,\d}$ defined by
\begin{align}
&T_{\g,\d}f = T_{max}f,\lb{3.2}\\
&f\in \dom(T_{\g,\d}) = \bigg\{g\in \dom(T_{max})\,\bigg|\,
\begin{aligned}
\cos(\g)\wti g(a)+\sin(\g)\wti g^{\, \prime}(a)&=0,\no\\
\cos(\d)\wti g(b)+\sin(\d)\wti g^{\, \prime}(b)&=0
\end{aligned}
\bigg\},\no
\end{align}
is a self-adjoint extension of $T_{min}$.\\[1mm]
$(ii)$  If $\varphi\in [0,\pi)$ and $R\in \SL(2,\bbR)$, then the operator $T_{\varphi,R}$ defined by
\begin{align}
&T_{\varphi,R}f = T_{max}f,\\
&f\in \dom(T_{\varphi,R}) = \bigg\{g\in \dom(T_{max})\,\bigg|\,
\begin{pmatrix}
\wti g(b)\\
\wti g^{\, \prime}(b)
\end{pmatrix}
= e^{i\varphi}R
\begin{pmatrix}
\wti g(a)\\
\wti g^{\, \prime}(a)
\end{pmatrix}
\bigg\},\no
\end{align}
is a self-adjoint extension of $T_{min}$.\\[1mm]
$(iii)$  If $T$ is a self-adjoint extension of $T_{min}$, then either $T=T_{\g,\d}$ for some $\g,\d\in[0,\pi)$, or $T=T_{\varphi,R}$ for some $\varphi\in [0,\pi)$ and some $R\in \SL(2,\bbR)$. \\[2mm] 
\noindent 
{\bf Notational Convention.} {\it To describe all possible self-adjoint boundary conditions associated with self-adjoint extensions of $T_{min}$ effectively, we will frequently employ the notation $T_{A,B}$, $M_{A,B}^{Do}(\dott)$, etc., where $A,B$ represents $\g,\d$ in the case of separated boundary conditions and $\varphi,R$ in the context of coupled boundary conditions.}

\smallskip

Choosing $\g=\d=0$ in \eqref{3.2} yields the self-adjoint extension with Dirichlet-type boundary conditions at $a$ and $b$, equivalently, the Friedrichs extension $T_F$ of $T_{min}$:
\begin{equation}
\dom(T_{0,0}) = \dom(T_F) = \{g\in\dom(T_{max})\, |\, \wti g(a) = \wti g(b)=0\}.
\end{equation}

Since the coefficients of the Sturm--Liouville differential expression are real, the following conjugation property holds:
\begin{equation}
\ol{u_j(z,\dott)} = u_j(\ol{z},\dott),\quad z\in \rho(T_{0,0}),\, j\in\{1,2\}.
\end{equation}
Applying \eqref{3.1}, one computes
\begin{equation}
\begin{split}
W(u_1(z,\dott),u_2(z,\dott)(a) &= -\wti u_1^{\, \prime}(z,a),\\
W(u_1(z,\dott),u_2(z,\dott)(b) &= \wti u_2^{\, \prime}(z,b),\quad z\in \rho(T_{0,0}).
\end{split}
\end{equation}
In particular, since the Wronskian of two solutions is constant,
\begin{equation}\lb{3.7}
\wti u_2^{\, \prime}(z,b) = -\wti u_1^{\, \prime}(z,a),\quad z\in \rho(T_{0,0}).
\end{equation}

We begin by recalling the orthonormal basis for $\cN_{\pm i}$ given by $\{v_j(\pm i,\dott)\}_{j=1,2}$,
\begin{align}
v_1(\pm i,\dott)& =c_1(\pm i)u_1(\pm i,\dott),\lb{3.8}\\
v_2(\pm i,\dott)& = c_2(\pm i)\bigg[u_2(\pm i,\dott)-\frac{(u_1(\pm i,\dott),u_2(\pm i,\dott))_{L^2((a,b);rdx)}}{\|u_1(\pm i,\dott)\|_{L^2((a,b);rdx)}^{2}}u_1(\pm i,\dott)\bigg]\\
&\;= c_2(\pm i)\bigg[u_2(\pm i,\dott)-\frac{\Im\big(\wti u_2^{\, \prime}(i,b)\big)}{\Im\big(\wti u_1^{\, \prime}(i,b) \big)}u_1(\pm i,\dott)\bigg],\no
\end{align}
with
\begin{align}
c_1(\pm i)& = \|u_1(\pm i,\dott)\|_{L^2((a,b);rdx)}^{-1}=\big[\mp \Im\big(\wti u_1^{\, \prime}(\pm i,b) \big)\big]^{-1/2},\\
c_2(\pm i)& = \bigg\|u_2(\pm i,\dott)-\frac{\Im\big(\wti u_2^{\, \prime}(i,b)\big)}{\Im\big(\wti u_1^{\, \prime}(i,b) \big)}u_1(\pm i,\dott)\bigg\|_{L^2((a,b);rdx)}^{-1}\lb{3.11}\\
&\;= \bigg[\pm \Im\big(\wti u_2^{\, \prime}(\pm i,a) \big)\pm\frac{\big[\Im\big(\wti u_2^{\, \prime}(\pm i,b)\big)\big]^2}{\Im\big(\wti u_1^{\, \prime}(\pm i,b) \big)}  \bigg]^{-1/2}. \no
\end{align}

The Donoghue $m$-function $M_{T_{A,B},\, \cN_i}^{Do}(\dott)$ with $T_{A,B}$ any self-adjoint extension of $T_{min}$ is provided next (cf. Theorems 6.1--6.3 in \cite{GLNPS21}).

\begin{theorem}\lb{t4.2}
Assume Hypothesis \ref{h4.1} and let $\{v_j(i,\dott)\}_{j=1,2}$ be the orthonormal basis for $\cN_i$ defined in \eqref{3.8}--\eqref{3.11}.  The Donoghue $m$-function $M_{T_{0,0},\, \cN_i}^{Do}(\dott):\bbC\backslash\bbR\to \cB(\cN_i)$ for $T_{0,0}$ satisfies 
\begin{align}
M_{T_{0,0},\, \cN_i}^{Do}(\pm i) &= \pm i I_{\cN_i},      \no \\
M_{T_{0,0},\, \cN_i}^{Do}(z) &= -\sum_{j,k=1}^2 [i \delta_{j,k} + W_{j,k}(z)] (v_k(i,\dott),\dott)_{L^2((a,b);rdx)}v_j(i,\dott)\big|_{\cN_i},      \\
&= -iI_{\cN_i}-\sum_{j,k=1}^2W_{j,k}(z)\big(v_k(i,\dott),\dott\big)_{L^2((a,b);rdx)}v_j(i,\dott)\big|_{\cN_i},  \no  \\
& \hspace*{6.2cm} z\in \bbC\backslash\bbR,\, z\neq \pm i,\no
\end{align}
where the matrix $\big(W_{j,k}(\dott)\big)_{j,k=1}^2$, $z\in \bbC\backslash\bbR$, $z\neq \pm i$, is given by
\begin{align}
&W_{1,1}(z)= [c_1(i)]^2\big[\wti u_1^{\, \prime}(z,b)-\wti u_1^{\, \prime}(-i,b)\big],\\
&W_{1,2}(z)= c_1(i)c_2(i)\Bigg\{\frac{\Im\big(\wti u_2^{\, \prime}( i,b)\big)}{\Im\big(\wti u_1^{\, \prime}( i,b) \big)}\big[\wti u_1^{\, \prime}(-i,b)-\wti u_1^{\, \prime}(z,b)\big]\\
&\hspace*{3.3cm}+\wti u_2^{\, \prime}(z,b)+\wti u_1^{\, \prime}(-i,a)\Bigg\}, \no\\
&W_{2,1}(z)= -c_1(i)c_2(i)\Bigg\{\frac{\Im\big(\wti u_2^{\, \prime}( i,b)\big)}{\Im\big(\wti u_1^{\, \prime}(i,b) \big)}\big[\wti u_1^{\, \prime}(z,b) - \wti u_1^{\, \prime}(-i,b)\big]\\
&\hspace*{3.5cm} + \wti u_2^{\, \prime}(-i,b)+\wti u_1^{\, \prime}(z,a)\Bigg\},\no\\
&W_{2,2}(z) = [c_2(i)]^2\Bigg\{\Bigg[\wti u_2^{\, \prime}(-i,b)-\wti u_2^{\, \prime}(z,b)\\
&\hspace*{3.1cm}+\frac{\Im\big(\wti u_2^{\, \prime}(i,b)\big)}{\Im\big(\wti u_1^{\, \prime}(i,b) \big)}\big[\wti u_1^{\, \prime}(z,b)-\wti u_1^{\, \prime}(-i,b)\big]\Bigg]\frac{\Im\big(\wti u_2^{\, \prime}(i,b)\big)}{\Im\big(\wti u_1^{\, \prime}(i,b) \big)}\no\\
&\hspace*{3cm} +\wti u_2^{\, \prime}(-i,a) - \wti u_2^{\, \prime}(z,a)+\frac{\Im\big(\wti u_2^{\, \prime}(i,b)\big)}{\Im\big(\wti u_1^{\, \prime}(i,b) \big)}\big[\wti u_1^{\, \prime}(z,a)-\wti u_1^{\, \prime}(-i,a)\big]\Bigg\}.\no
\end{align}
Furthermore, the following items $(i)$--$(v)$ hold.\\[2mm]
$(i)$ If $\g,\d\in (0,\pi)$, then the Donoghue $m$-function $M_{T_{\g,\d},\, \cN_i}^{Do}(\dott):\bbC\backslash\bbR\to \cB(\cN_i)$ for $T_{\g,\d}$ satisfies 
\begin{align}
& M_{T_{\g,\d},\, \cN_i}^{Do}(\pm i) = \pm i I_{\cN_i},    \no \\ 
& M_{T_{\g,\d},\, \cN_i}^{Do}(z)= M_{T_{0,0},\, \cN_i}^{Do}(z)\\
& \quad + (i-z)\sum_{j,k,\ell=1}^2 \big[K_{\g,\d}(z)^{-1}\big]_{j,k}W_{\ell,k}^{Kr}(z) ( u_j(\ol{z},\dott),\,\cdot\,)_{L^2((a,b);rdx)} v_{\ell}(i,\dott)\big|_{\cN_i},\no\\
&\hspace*{8.75cm} z\in\bbC\backslash\bbR,\, z\neq\pm i,\no
\end{align}
where the invertible matrix $K_{\g,\d}(\dott)$ and $\big(W_{\ell,k}^{Kr}(\dott)\big)_{\ell,k=1}^2$ are given by
\begin{align}
K_{\g,\d}(z) &=
\begin{pmatrix}
\cot(\d)+\wti u_1^{\, \prime}(z,b) & -\wti u_1^{\, \prime}(z,a)\\[2mm]
\wti u_2^{\, \prime}(z,b) & -\cot(\g)-\wti u_2^{\, \prime}(z,a)
\end{pmatrix},  \\
W_{1,1}^{Kr}(z)&=c_1(i)\big[\wti u_1^{\, \prime}(z,b) - \wti u_1^{\, \prime}(-i,b) \big],\lb{3.19}\\
W_{1,2}^{Kr}(z)&=c_1(i)\big[\wti u_2^{\, \prime}(z,b) + \wti u_1^{\, \prime}(-i,a) \big],\lb{3.20}\\
W_{2,1}^{Kr}(z)&=\wti v_2(-i,b)\wti u_1^{\, \prime}(z,b) - \wti v_2^{\, \prime}(-i,b) - \wti v_2(-i,a)\wti u_1^{\,\prime}(z,a)\lb{3.21}\\
&=-c_2(i)\bigg\{\frac{\Im\big(\wti u_2^{\, \prime}(i,b)\big)}{\Im\big(\wti u_1^{\, \prime}(i,b) \big)}\big[\wti u_1^{\, \prime}(z,b) - \wti u_1^{\, \prime}(-i,b)\big]+\wti u_2^{\, \prime}(-i,b)+\wti u_1^{\,\prime}(z,a)\bigg\},\no\\
W_{2,2}^{Kr}(z)&=\wti v_2(-i,b)\wti u_2^{\, \prime}(z,b) - \wti v_2(-i,a)\wti u_2^{\, \prime}(z,a) + \wti v_2^{\, \prime}(-i,a)\lb{3.22}\\
&=-c_2(i)\bigg\{\frac{\Im\big(\wti u_2^{\, \prime}(i,b)\big)}{\Im\big(\wti u_1^{\, \prime}(i,b) \big)}\big[\wti u_2^{\, \prime}(z,b) + \wti u_1^{\, \prime}(-i,a)\big] + \wti u_2^{\, \prime}(z,a) - \wti u_2^{\, \prime}(-i,a)\bigg\}.\no
\end{align}
$(ii)$ If $\varphi\in[0,\pi)$ and $R\in \SL(2,\bbR)$ with $R_{1,2}\neq 0$, then the Donoghue $m$-function $M_{T_{\varphi,R},\, \cN_i}^{Do}(\dott):\bbC\backslash\bbR\to \cB(\cN_i)$ for $T_{\varphi,R}$ satisfies
\begin{align}
& M_{T_{\varphi,R},\, \cN_i}^{Do}(\pm i) = \pm i I_{\cN_i},    \no \\ 
& M_{T_{\varphi, R},\, \cN_i}^{Do}(z)= M_{T_{0,0},\, \cN_i}^{Do}(z)\\
& \quad + (i-z)\sum_{j,k,\ell=1}^2 \big[K_{\varphi,R}(z)^{-1}\big]_{j,k}W_{\ell,k}^{Kr}(z) ( u_j(\ol{z},\dott),\,\cdot\,)_{L^2((a,b);rdx)} v_{\ell}(i,\dott)\big|_{\cN_i},\no\\
&\hspace*{8.85cm} z\in\bbC\backslash\bbR,\, z\neq\pm i,\no
\end{align}
where $\big(W_{\ell,k}^{Kr}(\dott)\big)_{\ell,k=1}^2$ is once again given in \eqref{3.19}--\eqref{3.22} and the invertible matrix $K_{\varphi,R}(\dott)$ is given by
\begin{equation}
K_{\varphi,R}(z) = \begin{pmatrix}
-\dfrac{R_{2,2}}{R_{1,2}}+\wti u_1^{\, \prime}(z,b) & \dfrac{e^{-i\varphi}}{R_{1,2}}-\wti u_1^{\, \prime}(z,a)\\[4mm]
\dfrac{e^{i\varphi}}{R_{1,2}}+\wti u_2^{\, \prime}(z,b) & -\dfrac{R_{1,1}}{R_{1,2}}-\wti u_2^{\, \prime}(z,a)
\end{pmatrix}.
\end{equation}
$(iii)$ If $\g\in (0,\pi)$, then the Donoghue $m$-function $M_{T_{\g,0},\, \cN_i}^{Do}(\dott):\bbC\backslash\bbR\to \cB(\cN_i)$ for $T_{\g,0}$ satisfies 
\begin{align}
M_{T_{\g,0},\, \cN_i}^{Do}(\pm i) &= \pm i I_{\cN_i},    \no \\
M_{T_{\g,0},\, \cN_i}^{Do}(z) &= M_{T_{0,0},\, \cN_i}^{Do}(z)\\
&\quad + \frac{z-i}{\cot(\g)+\wti u_2^{\, \prime}(z,a)} ( u_2(\ol{z},\dott),\,\cdot\,)_{L^2((a,b);rdx)} \sum_{\ell=1}^2 W_{\ell,2}^{Kr}(z)v_{\ell}(i,\dott)\big|_{\cN_i},\no\\
&\hspace*{7.75cm} z\in\bbC\backslash\bbR,\, z\neq\pm i,\no
\end{align}
where $\cot(\g)+\wti u_2^{\, \prime}(z,a)\neq0$ and the scalars $\big\{W_{\ell,2}^{Kr}(\dott)\big\}_{\ell=1,2}$ are given by \eqref{3.20} and \eqref{3.22}.\\[2mm]
$(iv)$ If $\d\in (0,\pi)$, then the Donoghue $m$-function $M_{T_{0,\d},\, \cN_i}^{Do}(\dott):\bbC\backslash\bbR\to \cB(\cN_i)$ for $T_{0,\d}$ satisfies 
\begin{align}
M_{T_{0,\d},\, \cN_i}^{Do}(\pm i) &= \pm i I_{\cN_i},     \no \\
M_{T_{0,\d},\, \cN_i}^{Do}(z) &= M_{T_{0,0},\, \cN_i}^{Do}(z)\\
&\quad - \frac{z-i}{\cot(\d)+\wti u_1^{\, \prime}(z,b)} ( u_1(\ol{z},\dott),\,\cdot\,)_{L^2((a,b);rdx)} \sum_{\ell=1}^2 W_{\ell,1}^{Kr}(z)v_{\ell}(i,\dott)\big|_{\cN_i},\no\\
&\hspace*{7.65cm} z\in\bbC\backslash\bbR,\, z\neq\pm i,\no
\end{align}
where $\cot(\d)+\wti u_1^{\, \prime}(z,b)\neq0$ and the scalars $\big\{W_{\ell,1}^{Kr}(\dott)\big\}_{\ell=1,2}$ are given by \eqref{3.19} and \eqref{3.21}.\\[2mm]
$(v)$  If $\varphi \in [0,\pi)$ and $R\in \SL(2,\bbR)$ with $R_{1,2}=0$, then the Donoghue $m$-function $M_{T_{\varphi,R},\, \cN_i}^{Do}(\dott):\bbC\backslash\bbR\to \cB(\cN_i)$ for $T_{\varphi,R}$ satisfies 
\begin{align}
& M_{T_{\varphi,R},\, \cN_i}^{Do}(\pm i)=\pm i I_{\cN_i},     \no \\
&M_{T_{\varphi,R},\, \cN_i}^{Do}(z) = M_{T_{0,0},\, \cN_i}^{Do}(z)      \\
&\quad- \frac{z-i}{k_{\varphi,R}(z)} (u_{\varphi,R}(\ol{z},\dott),\dott)_{L^2((a,b);rdx)}\sum_{\ell=1}^2\big[e^{-i\varphi}R_{2,2}W_{\ell,2}^{Kr}(z)+W_{\ell,1}^{Kr}(z) \big]v_{\ell}(i,\dott)\big|_{\cN_i},\no\\
&\hspace*{10cm} z\in\bbC\backslash\bbR,\, z\neq\pm i,\no
\end{align}
where the matrix $\big(W_{\ell,k}^{Kr}(\dott)\big)_{\ell,k=1}^2$ is once again given in \eqref{3.19}--\eqref{3.22} and the nonzero scalar $k_{\varphi,R}(\dott)$ is given by
\begin{align}
k_{\varphi,R}(z)=-R_{2,1}R_{2,2}-e^{i\varphi}R_{2,2}\wti u_{\varphi,R}^{\, \prime}(z,a)+\wti u_{\varphi,R}^{\, \prime}(z,b),
\end{align}
where
\begin{align}
u_{\varphi,R}(\zeta,\dott) = e^{-i\varphi}R_{2,2}u_2(\zeta,\dott)+u_1(\zeta,\dott),\quad \zeta\in \rho(T_{0,0}).
\end{align}
\end{theorem}

\begin{remark}
For the Krein extension, $T_{0,R_K}$, under the additional assumption that $T_{min} \geq \varepsilon I_{_{L^2((a,b);rdx)}}$ for some 
$\varepsilon > 0$, applying \cite[Theorem~3.5\,$(ii)$]{FGKLNS21}, one computes for the matrix $K_{0,R_K}$,
\begin{equation}\lb{3.25}
K_{0,R_K}(z) = \begin{pmatrix}
\wti u_1^{\, \prime}(z,b)-\wti u_1^{\, \prime}(0,b) & \wti u_1^{\, \prime}(0,a)-\wti u_1^{\, \prime}(z,a)\\[4mm]
\wti u_2^{\, \prime}(z,b)-\wti u_2^{\, \prime}(0,b) & \wti u_2^{\, \prime}(0,a)-\wti u_2^{\, \prime}(z,a)
\end{pmatrix}, \quad z\in \rho(T_{0,0})\cap\rho(T_{0,R_K}), 
\end{equation}
where we note that $0 \in \sigma(T_{0,R_K})$. \hfill$\diamond$
\end{remark}

\section{Donoghue $m$-functions: One Limit Circle Endpoint} \lb{s3}

In this section we recall the Donoghue $m$-functions in the case where $\tau$ is in the limit circle case at precisely one endpoint (which we choose to be $a$ without loss of generality) following \cite[Sect. 5]{GLNPS21}. 

\begin{hypothesis}\lb{h3.1}
In addition to Hypothesis \ref{h2.1} assume that $\tau$ is in the limit circle case at $a$ and in the limit point case at $b$. 
Moreover, for $z\in \rho(T_0)$, let $\psi(z,\dott)$ denote the unique solution to $(\tau -z)y=0$ that satisfies $\psi(z,\dott)\in L^2((a,b);rdx)$ and $\wti \psi(z,a)=1$.
\end{hypothesis}

Assume Hypothesis \ref{h3.1}. By Theorem \ref{t2.8} or Theorem \ref{t2.13}, the following statements $(i)$ and $(ii)$ hold.\\[1mm]
$(i)$  If $\g\in [0,\pi)$, then the operator $T_{\g}$ defined by
\begin{equation}\lb{4.1}
\begin{split}
&T_{\g}f=T_{max}f,\\
&f\in \dom(T_{\g})=\{g\in \dom(T_{max})\,|\, \cos(\g)\wti g(a) + \sin(\g)\wti g^{\, \prime}(a)=0\}, 
\end{split}
\end{equation}
is a self-adjoint extension of $T_{min}$.\\[1mm]
$(ii)$  If $T$ is a self-adjoint extension of $T_{min}$, then $T=T_{\g}$ for some $\g\in[0,\pi)$.\\[2mm]
Statements analogous to $(i)$ and $(ii)$ hold if $\tau$ is in the limit point case at $a$ and in the limit circle case at $b$; for brevity we omit the details.

Choosing $\g=0$ in \eqref{4.1} yields the self-adjoint extension $T_0$ with a Dirichlet-type boundary condition at $a$:
\begin{equation}
\dom(T_0)=\{g\in\dom(T_{max})\, |\, \wti g(a) = 0\}.
\end{equation}

Since the coefficients $p$, $q$, and $r$ are real-valued, the solution $\psi(z,\dott)$ has the following conjugation property:
\begin{equation}
\overline{\psi(z,\dott)} = \psi(\overline{z},\dott),\quad z\in \rho(T_0).
\end{equation}

We now turn to the Donoghue $m$-function $M_{T_\g,\, \cN_i}^{Do}(\dott)$ with $T_\g$ any self-adjoint extension of $T_{min}$ (cf. Theorems 5.1 and 5.2 in \cite{GLNPS21}).

\begin{theorem}\lb{t3.2}
Assume Hypothesis \ref{h3.1} and let $\g\in [0,\pi)$.  The Donoghue $m$-function $M_{T_\g,\, \cN_i}^{Do}(\dott):\bbC\backslash\bbR\to \cB(\cN_i)$
for $T_{\g}$ satisfies 
\begin{align}
M_{T_\g,\, \cN_i}^{Do}(\pm i) &= \pm i I_{\cN_i}, \quad \g\in[0,\pi),   \no \\
M_{T_0,\, \cN_i}^{Do}(z) &= \Bigg[-i +  \frac{\wti \psi^{\, \prime}(z,a)-\wti \psi^{\, \prime}(-i,a)}{\Im\big(\wti \psi^{\, \prime}(i,a)\big)}\Bigg]I_{\cN_i},\quad z\in \bbC\backslash\bbR,\, z\neq \pm i,  \no  \\
M_{T_\g,\, \cN_i}^{Do}(z) &= M_{T_0,\, \cN_i}^{Do}(z) \\
&\quad + (i-z)\frac{\wti \psi^{\, \prime}(z,a)-\wti \psi^{\, \prime}(-i,a)}{\cot(\g) + \wti \psi^{\, \prime}(z,a)} (\psi(\overline{z},\dott),\dott)_{L^2((a,b);rdx)}\psi(i,\dott)\bigg|_{\cN_i},\no\\
&\hspace*{5.55cm} \g\in(0,\pi),\, z\in\bbC\backslash\bbR,\, z\neq \pm i.  \no
\end{align}
\end{theorem}

\section{The Jacobi Operator and its Donoghue $m$-functions} \lb{s6}

We now turn to the principal topic of this paper, the Jacobi differential expression 
\begin{align} \lb{6.1}
\begin{split} 
\tau_{\a,\b} = - (1-x)^{-\a} (1+x)^{-\b}(d/dx) \big((1-x)^{\a+1}(1+x)^{\b+1}\big) (d/dx),&     \\ 
x \in (-1,1), \; \a, \b \in \bbR, 
\end{split} 
\end{align}
that is, in connection with Sections \ref{s2} one now has  
\begin{align}
& a=-1, \quad b = 1,    \no \\
& p(x) = p_{\a,\b}(x) = (1-x)^{\a+1}(1+x)^{\b+1}, \quad q (x) = q_{\a,\b}(x) = 0,  \\ 
& r(x) = r_{\a,\b}(x) = (1-x)^\a (1+x)^\b, \quad x\in(-1,1),\quad \a,\b\in\bbR     \no 
\end{align}
(see, e.g. \cite[Ch.~22]{AS72}, \cite{BEZ01}, \cite{Bo29}, \cite[Sect.~23]{Ev05}, \cite[Ch.~4]{Is05}, \cite[Sects.~VII.6.1, XIV.2]{Kr02}, \cite[Ch.~18]{Ov20}, \cite[Ch.~IV]{Sz75}). 

$L^2$-realizations of $\tau_{\a,\b}$ are thus most naturally associated with the Hilbert space 
$L^2((-1,1); r_{\a,\b} dx)$. However, occasionally the weight function is absorbed into the Hilbert space leading to an equivalent differential expression in the Hilbert space $L^2((-1,1);dx)$ (cf. \cite[p.~1510--1520]{DS88}, \cite[Sect.~37]{Ev05}, \cite{Gr74}). For more recent developments see, for instance, \cite{EKLWY07}, \cite{Fr20}, \cite{FL20}, \cite{KKT18}, \cite{KMO05}.

To decide the limit point/limit circle classification of $\tau_{\a,\b}$ at the interval endpoints $\pm 1$, it suffices 
to note that if $y_1$ is a given solution of $\tau y = 0$, then a 2nd linearly independent solution $y_2$ of 
$\tau y = 0$ is obtained via the standard formula
\begin{equation}
y_2(x) = y_1(x) \int_c^x dx' \, p(x')^{-1} y_1(x')^{-2}, \quad c, x \in (a,b).    
\end{equation}
Returning to the concrete Jacobi case at hand, one can choose
\begin{align} 
& y_1(x) = 1, \quad x \in (-1,1), \no \\
& y_2(x) = \int_0^x dx' \, (1 - x')^{-1-\a} (1+x')^{-1-\b}      \lb{6.5} \\
&\quad = \begin{cases}
2^{-1-\a} \b^{-1} (1+ x)^{-\b} [1+\Oh(1+ x)] + \Oh(1), & \a \in \bbR, \, \b \in \bbR\backslash\{0\}, \; \text{as $x \downarrow -1$}, \\
- 2^{-1-\a} \ln(1+ x) + \Oh(1), & \a \in \bbR, \, \b = 0, \; \text{as $x \downarrow -1$},  \\
2^{-1-\b} \a^{-1} (1 - x)^{-\a} [1+\Oh(1- x)] + \Oh(1), & \a \in \bbR\backslash\{0\}, \, \b \in \bbR, \; \text{as $x \uparrow +1$}, \\
- 2^{-1-\b} \ln(1 - x) + \Oh(1), & \a =0, \, \b \in \bbR, \; \text{as $x \uparrow +1$}.   
\end{cases}    \no 
\end{align} 
Thus, one has the classification,
\begin{equation}
\tau_{\a,\b} \, \text{ is } \begin{cases} \text{regular at $-1$ if and only if $\a \in \bbR$, $\b \in (-1,0)$,} \\
\text{in the limit circle case at $-1$ if and only if $\a \in \bbR$, $\b \in [0,1)$,} \\
\text{in the limit point case at $-1$ if and only if $\a \in \bbR$, $\b \in \bbR \backslash (-1,1)$,} \\
\text{regular at $+1$ if and only if $\a \in (-1,0)$, $\b \in \bbR$,} \\
\text{in the limit circle case at $+1$ if and only if $\a \in [0,1)$, $\b \in \bbR$,} \\
\text{in the limit point case at $+1$ if and only if $\a \in \bbR \backslash (-1,1)$, $\b \in \bbR$.}
\end{cases}     \lb{6.6}
\end{equation}
The maximal and preminimal operators, $T_{max,\a,\b}$ and $T_{min,0,\a,\b}$, associated to $\tau_{\a,\b}$ 
in $L^2((-1,1); r_{\a,\b} dx)$ are then given by
\begin{align}
&T_{max,\a,\b} f = \tau_{\a,\b} f,     \no
\\
& f \in \dom(T_{max,\a,\b})=\big\{g\in L^2((-1,1); r_{\a,\b} dx) \, \big| \,  g,g^{[1]}\in AC_{loc}((-1,1)); \no \\
& \hspace*{6.8cm} \tau_{\a,\b} g\in L^2((-1,1); r_{\a,\b} dx)\big\},
\end{align}
and
\begin{align}
&T_{min,0,\a,\b} f = \tau_{\a,\b} f, \no
\\
& f \in \dom(T_{min,0,\a,\b})=\big\{g\in L^2((-1,1); r_{\a,\b} dx)  \, \big| \,  g,g^{[1]}\in AC_{loc}((-1,1));    \\ 
&\hspace*{2.2cm} \supp \, (g)\subset(-1,1) \text{ is compact; } 
\tau_{\a,\b} g\in L^2((-1,1); r_{\a,\b} dx)\big\}.   \no 
\end{align}
The fact \eqref{6.5} naturally leads to principal and nonprincipal solutions 
$u_{\pm1,\a,\b}(0,x)$ and $\hatt u_{\pm1,\a,\b}(0,x)$ of $\tau_{\a,\b}y=0$ near $\pm1$ as follows: 
\begin{align} 
\begin{split} 
u_{-1,\a,\b}(0,x)&=\begin{cases}
-2^{-\a-1}\b^{-1} (1+x)^{-\b}[1+\Oh(1+x)], & \b\in(-\infty,0),\\
1, & \b\in[0,\infty),
\end{cases}     \\
\hatt u_{-1,\a,\b}(0,x)&=\begin{cases}
1, & \b\in(-\infty,0),\\
-2^{-\a-1}\ln((1+x)/2), & \b=0,\\
2^{-\a-1} \b^{-1} (1+x)^{-\b} [1+\Oh(1+x)], & \b\in(0,\infty),
\end{cases}
\end{split} 
\quad \a\in\bbR,     \lb{6.9} 
\end{align}
and
\begin{align}
\begin{split} 
u_{+1,\a,\b}(0,x)&=\begin{cases}
2^{-\b-1} \a^{-1} (1-x)^{-\a} [1+\Oh(1-x)], & \a\in(-\infty,0),\\
1, & \a\in[0,\infty),
\end{cases}      \\
\hatt u_{+1,\a,\b}(0,x)&=\begin{cases}
1, & \a\in(-\infty,0),\\
2^{-\b-1}\ln((1-x)/2), & \a=0,\\
-2^{-\b-1} \a^{-1} (1-x)^{-\a} [1+\Oh(1-x)], & \a\in(0,\infty),
\end{cases}
\end{split} 
\quad \b\in\bbR.     \lb{6.10} 
\end{align}

Combining the fact \eqref{6.6} with Theorem \ref{t2.5}, $T_{min,0,\a,\b}$ is essentially self-adjoint in 
$L^2((-1,1); r_{\a,\b} dx)$ if and only if 
$\a, \b \in \bbR \backslash (-1,1)$. Thus, boundary values for $T_{max,\a,\b}$ at $-1$ exist if and only if 
$\a \in \bbR$, $\b \in (-1,1)$, and similarly, boundary values for $T_{max,\a,\b}$ at $+1$ exist if and only if 
$\a \in (-1,1)$, $\b \in \bbR$. 

Employing the principal and nonprincipal solutions \eqref{6.9}, \eqref{6.10} at $\pm 1$, according to 
\eqref{2.25}, \eqref{2.26}, generalized boundary values for $g\in\dom(T_{max,\a,\b})$ are of the form
\begin{align}
\begin{split} 
\wti g(-1)&=\begin{cases}
g(-1), & \b\in(-1,0),\\
-2^{\a+1}\lim_{x\downarrow-1}g(x)/\ln((1+x)/2), & \b=0,\\
\b 2^{\a+1}\lim_{x\downarrow-1}(1+x)^\b g(x), & \b\in(0,1),
\end{cases}    \\
{\wti g}^{\, \prime}(-1)&=\begin{cases}
g^{[1]}(-1), &\b\in(-1,0),\\
\lim_{x\downarrow-1}\big[g(x)+\wti g(-1)2^{-\a-1}\ln((1+x)/2)\big], & \b=0,\\
\lim_{x\downarrow-1}\big[g(x)-\wti g(-1)2^{-\a-1}\b^{-1}(1+x)^{-\b}\big], & \b\in(0,1),
\end{cases}
\end{split} 
\quad \a\in\bbR,     \\
\begin{split}
\wti g(1)&=\begin{cases}
g(1), & \a\in(-1,0),\\
2^{\b+1}\lim_{x\uparrow1}g(x)/\ln((1-x)/2), & \a=0,\\
-\a2^{\b+1}\lim_{x\uparrow1}(1-x)^\a g(x), & \a\in(0,1),
\end{cases}\\
{\wti g}^{\, \prime}(1)&=\begin{cases}
g^{[1]}(1), & \a\in(-1,0),\\
\lim_{x\uparrow1}\big[g(x)-\wti g(1)2^{-\b-1}\ln((1-x)/2)\big], & \a=0,\\
\lim_{x\uparrow1}\big[g(x)+\wti g(1)2^{-\b-1}\a^{-1}(1-x)^{-\a}\big], & \a\in(0,1),
\end{cases}
\end{split}
\quad \b\in\bbR.
\end{align}

As a result, the minimal operator 
$T_{min}$ associated to $\tau_{\a,\b}$, that is, $T_{min} = \ol{T_{min,0}}$, is thus given by 
\begin{align}
&T_{min,\a,\b} f = \tau_{\a,\b} f, \no \\
& f \in \dom(T_{min,\a,\b})=\big\{g\in L^2((-1,1); r_{\a,\b} dx) \, \big| \, g,g^{[1]}\in AC_{loc}((-1,1));    \\ 
&\hspace*{1.1cm} \wti g(-1) = {\wti g}^{\, \prime}(-1) = \wti g(1) = {\wti g}^{\, \prime}(1) = 0; \, 
\tau_{\a,\b} g\in L^2((-1,1); r_{\a,\b} dx)\big\}.  \no
\end{align}

For a detailed treatment of solutions of the Jacobi differential equation and the associated hypergeometric differential equations we refer to Appendices \ref{sA}--\ref{sC}.

\begin{remark}
We now mention a few special cases of interest. The Legendre equation ($\a=\b=0$) has frequently been discussed in the literature, 
see, for instance, \cite{GLN20} and the extensive list of references cited therein. The Gegenbauer, or ultraspherical, equation (see, e.g., \cite[Ch. 22]{AS72}, \cite[Ch.~18]{Ov20}, \cite[Ch. IV]{Sz75}) can be realized by choosing the parameters $\a=\b=\mu-1/2$, noting  at the endpoints $x=\pm1$, $\tau_\mu$ is regular for $\mu\in(-1/2,1/2)$, in the limit circle case for $\mu\in[1/2,3/2)$, and in the limit point case for $\mu\in\bbR\backslash(-1/2,3/2)$. The Chebyshev equations of the first and second kinds are two more important special cases, with the first kind realized by choosing $\mu=0$ in the Gegenbauer equation, or $\a=\b=-1/2$ in the Jacobi equation (see, e.g., \cite[Ch. 22]{AS72}, \cite[Ch.~18]{Ov20}, \cite[Ch. IV]{Sz75}), whereas the second kind is realized by choosing $\mu=1$ in the Gegenbauer equation, or $\a=\b=1/2$ in the Jacobi equation (see, e.g., \cite[Ch. 22]{AS72}, \cite[Ch.~18]{Ov20}, \cite[Ch. IV]{Sz75}).
\hfill$\diamond$
\end{remark}

We now determine the solutions $\phi_{0,\a,\b}(z,\dott)$ and $\theta_{0,\a,\b}(z,\dott)$ of $\tau_{\a,\b}u=zu,z\in\bbC$, that are subject to the conditions
\begin{align}
\begin{split}
\wti \phi_{0,\a,\b} (z,-1) &= 0, \quad \wti \phi_{0,\a,\b}^{\, \prime} (z,-1) = 1,\\
\wti \theta_{0,\a,\b} (z,-1) &= 1, \quad \wti \theta_{0,\a,\b}^{\, \prime} (z,-1) = 0.    \lb{6.14} 
\end{split}
\end{align}
In particular, one obtains from $\eqref{C.13}$,
\begin{align}
\begin{split}
\phi_{0,\a,\b} (z,x) &=\begin{cases}
-2^{-\a-1}\b^{-1} y_{2,\a,\b,-1}(z,x), & \b\in(-1,0),\\
y_{1,\a,\b,-1}(z,x), & \b\in[0,1),\\
\end{cases}\\
\theta_{0,\a,\b} (z,x) &=\begin{cases}
y_{1,\a,\b,-1}(z,x), & \b\in(-1,0),\\
-2^{-\a-1}y_{2,\a,0,-1}(z,x), & \b=0,\\
2^{-\a-1}\b^{-1}y_{2,\a,\b,-1}(z,x), & \b\in(0,1),
\end{cases}\\
&\hspace*{2.4cm} \a\in\bbR,\; z\in\bbC,\; x\in(-1,1).
\end{split}
\end{align}

\subsection{The Regular and Limit Circle Case $\a, \b \in (-1,1)$}

In this section we compute the Donoghue $m$-function when the Jacobi problem considered is either in the regular or limit circle case at $\pm 1$. 

Using \eqref{6.14}, the solutions in \eqref{3.1} for this example are given by 
\begin{align}
&u_{1,\a,\b}(z,x)=\phi_{0,\a,\b}(z, x)/\wti\phi_{0,\a,\b}(z, 1) \no \\
&\hspace*{1.71cm} = \begin{cases}
y_{2,\a,\b,-1}(z,x)/\wti y_{2,\a,\b,-1}(z,1), & \b\in(-1,0),\\
y_{1,\a,\b,-1}(z,x)/\wti y_{1,\a,\b,-1}(z,1), & \b\in[0,1),
\end{cases}   \lb{6.16}  \\
&u_{2,\a,\b}(z,x)=\theta_{0,\a,\b}(z, x)-[\wti\theta_{0,\a,\b}(z, 1)/\wti\phi_{0,\a,\b}(z, 1)]\phi_{0,\a,\b}(z, x)   \no  \\
&=\begin{cases}
y_{1,\a,\b,-1}(z,x)-[\wti y_{1,\a,\b,-1}(z,1)/\wti y_{2,\a,\b,-1}(z,1)] y_{2,\a,\b,-1}(z,x), \\
\hspace*{7.5cm} \b\in(-1,0),\\
-2^{-\a-1}\{y_{2,\a,0,-1}(z,x)-[\wti y_{2,\a,0,-1}(z,1)/\wti y_{1,\a,0,-1}(z,1)] y_{1,\a,0,-1}(z,x)\},\\
\hspace*{9.75cm} \b=0,\\
2^{-\a-1}\b^{-1}\{y_{2,\a,\b,-1}(z,x)-[\wti y_{2,\a,\b,-1}(z,1)/\wti y_{1,\a,\b,-1}(z,1)] y_{1,\a,\b,-1}(z,x)\}, \\
\hfill \b\in(0,1),
\end{cases}  \no \\
&\hspace*{7.1cm} \a \in (-1,1),\; z\in\bbC,\; x \in (-1,1), \no
\end{align}
where the generalized boundary values are given in \eqref{C.14}--\eqref{C.16}.
Hence substituting \eqref{6.16} into \eqref{3.8}--\eqref{3.11} and applying Theorem \ref{t4.2} yields the (Nevanlinna--Herglotz) Donoghue $m$-function $M_{T_{A,B,\a,\b},\, \cN_i}^{Do}(\dott)$ for any self-adjoint extension $T_{A,B,\a,\b}$ of $T_{min}$ with $\a,\b\in(-1,1)$.

As an example of coupled boundary conditions, we consider the Krein--von Neumann extension following Example 4.3 found in \cite{FGKLNS21}. For $\a,\b\in(-1,1)$, the following five cases are associated with a strictly positive minimal operator 
$T_{min,\a,\b}$ and we now provide the corresponding choices of $R_{K,\a,\b}$ for the Krein--von Neumann extension $T_{0,R_K,\a,\b}$ of $T_{min,\a,\b}$:
\begin{align}
& T_{0,R_K,\a,\b} f = \tau_{\a,\b} f,    \\
& f \in \dom(T_{0,R_K,\a,\b})=\bigg\{g\in\dom(T_{max,\a,\b}) \, \bigg| \begin{pmatrix} \wti g(1) 
\\ {\wti g}^{\, \prime}(1) \end{pmatrix} = R_{K,\a,\b} \begin{pmatrix}
\wti g(-1) \\ {\wti g}^{\, \prime}(-1) \end{pmatrix} \bigg\}, \no \\
& R_{K,\a,\b}=\begin{cases}
\begin{pmatrix}  1 & 2^{-\a-\b-1}\dfrac{\Gamma(-\a)\Gamma(-\b)}{\Gamma(-\a-\b)} \\
0 & 1
\end{pmatrix}, & \a,\b\in(-1,0),\\[7mm]
\begin{pmatrix}  -2^{-\a-\b-1}\dfrac{\Gamma(-\a)\Gamma(-\b)}{\Gamma(-\a-\b)} & 1 \\
-1 & 0
\end{pmatrix}, & \a\in(-1,0),\; \b\in(0,1), \\[7mm]
\begin{pmatrix}  0 & -1 \\
1 & 2^{-\a-\b-1}\dfrac{\Gamma(-\a)\Gamma(-\b)}{\Gamma(-\a-\b)}
\end{pmatrix}, & \a\in(0,1),\; \b\in(-1,0), \\[7mm]
\begin{pmatrix}  0 & -1 \\
1 & -2^{-\b-1}[\gamma_{E}+\psi(-\b)]
\end{pmatrix}, & \a=0,\; \b\in(-1,0) \\[7mm]
\begin{pmatrix}  2^{-\a-1}[\gamma_{E}+\psi(-\a)] & 1 \\
-1 & 0
\end{pmatrix}, & \a\in(-1,0),\; \b=0,
\end{cases}\lb{6.18}
\end{align}
where we interpret $1/\Gamma(0) = 0$, $\psi(\dott) = \Gamma'(\dott)/\Gamma(\dott)$ denotes the Digamma function, and $\gamma_{E} = - \psi(1) = 0.57721\dots$ represents Euler's constant. Obviously, $\det(R_{K,\a,\b}) = 1$ in all five cases. Furthermore, as $R_{1,2}\neq0$ for each case, Theorem \ref{t4.2} $(ii)$ applies and one obtains the Donoghue $m$-function $M_{T_{0,R_K,\a,\b},\, \cN_i}^{Do}(\dott)$ for the Krein--von Neumann extension $T_{0,R_K,\a,\b}$ by utilizing \eqref{6.16} and \eqref{6.18} as well as the explicit form of $K_{0,R_K}(\dott)$ in \eqref{3.25}. We note once again that $M_{T_{0,R_K,\a,\b},\, \cN_i}^{Do}(\dott)$ is a Nevanlinna--Herglotz function.

In the remaining four cases not covered by \eqref{6.18}, given by all combinations of $\a=0,\ \b=0,\ \a\in(0,1)$, and 
$\b\in(0,1)$, one observes that \cite[Theorem 3.5]{FGKLNS21} is not applicable as the underlying minimal operator, $T_{min, \a,\b}$, is nonnegative but not strictly positive. In particular, the Jacobi polynomials satisfy Friedrichs boundary conditions for $\a,\b\in[0,1)$, hence $0\in\sigma(T_{F,\a,\b}),\ \a,\b\in[0,1)$ and $T_{min, \a,\b} \geq 0$ is nonnegative, but not strictly positive when 
$\a,\b \in [0,1)$.

\subsection{Precisely One Interval Endpoint in the Limit Point Case}

In this section we determine the Donoghue $m$-function in all situations where precisely one interval endpoint is in the limit point case. We will focus on the case when $\a \in (-\infty, -1]$ or 
$\a \in [1,\infty)$, so that the right endpoint $x=1$ represents the limit point case. The converse situation can be obtained by reflection with respect to the origin (i.e., considering the transform $(-1,1) \ni x \mapsto -x \in (-1,1)$).  

We recall from \cite[Sect.~6]{GLN20} that the Weyl--Titchmarsh--Kodaira solution and $m$-function corresponding to the Friedrichs (resp., Dirichlet) boundary condition at $x=-1$ is determined via the requirement
\begin{align}
\begin{split} 
\psi_{0,\a,\b} (z,\dott) = \theta_{0,\a,\b}(z,\dott) + m_{0,\a,\b} (z) \phi_{0,\a,\b}(z,\dott) 
\in L^2((c,1); r_{\a,\b} dx),&  \\
z \in \bbC \backslash \sigma(T_{F,\a,\b}), \;  \a \in (-\infty,-1]\cup[1,\infty), \; \b\in(-1,1), \; c \in (-1,1).&
\end{split}
\end{align}
In particular, since $\wti \psi_{0,\a,\b}^{\, \prime}(z,-1)=m_{0,\a,\b}(z)$ one finds from Theorem \ref{t3.2},
\begin{align}
\begin{split}
M_{T_{0,\a,\b},\, \cN_i}^{Do}(z)&=\Bigg[-i+\dfrac{m_{0,\a,\b}(z)-m_{0,\a,\b}(-i)}{\Im(m_{0,\a,\b}(i))}\Bigg]I_{\cN_i},  \\
M_{T_{\g,\a,\b},\, \cN_i}^{Do}(z) &= M_{T_{0,\a,\b},\, \cN_i}^{Do}(z)+ (i-z)\frac{m_{0,\a,\b}(z)-m_{0,\a,\b}(-i)}{\cot(\g) + m_{0,\a,\b}(z)}  \\
&\quad\; \times (\psi_{0,\a,\b}(\overline{z},\dott),\dott)_{L^2((a,b);rdx)}\psi_{0,\a,\b}(i,\dott)\big|_{\cN_i},\quad \g\in(0,\pi), \\
&\hspace{2.1cm}\a\in(-\infty,-1]\cup[1,\infty),\; \b\in(-1,1),\; z\in\bbC\backslash\bbR,
\end{split}
\end{align}
where $\psi_{0,\a,\b}(z,\dott)$ and $m_{0,\a,\b}(z,\dott)$ are given by the following:\\[1mm]
\noindent 
{\boldmath $\mathbf{(I)}$ {\bf The Case $\a \in [1,\infty)$ and $\b\in(-1,0)$:}}
\begin{align}
\begin{split}
& \hspace*{-.3cm} \psi_{0,\a,\b}(z,x)=y_{1,\a,\b,-1}(z,x)-2^{-\a-1}\b^{-1}y_{2,\a,\b,-1}(z,x) m_{0,\a,\b}(z),
\\
& m_{0,\a\,\b}(z) = 2^{1+\a+\b}\b \f{\Gamma(1+\b)}{\Gamma(1-\b)}     \\
& \hspace*{1.8cm}  \times \f{\Gamma([1+\a-\b+\sigma_{\a,\b}(z)]/2)
\Gamma([1+\a-\b-\sigma_{\a,\b}(z)]/2)}{
\Gamma([1+\a+\b+\sigma_{\a,\b}(z)]/2)\Gamma([1+\a+\b-\sigma_{\a,\b}(z)]/2)},     \\
&\hspace*{4.7cm}  z\in\rho(T_{F,\a,\b}), \; \a \in [1,\infty), \; \b\in(-1,0),    \\
&\sigma(T_{F,\a,\b}) = \{(n-\b)(n+1+\a)\}_{n\in\bbN_0}, \quad \a\in [1,\infty), \; \b\in(-1,0), 
\end{split}\lb{6.27}
\end{align}
with
\begin{equation}
\sigma_{\a,\b}(z) = \big[(1+\a+\b)^2+4z\big]^{1/2}.
\end{equation}
\noindent 
{\boldmath $\mathbf{(II)}$ {\bf The Case $\a \in [1,\infty)$ and $\b=0$:}} 
\begin{align}
& \hspace*{-.3cm} \psi_{0,\a,0}(z,x)=-2^{-\a-1}y_{2,\a,0,-1}(z,x)+y_{1,\a,0,-1}(z,x) m_{0,\a,0}(z), \no \\
& m_{0,\a,0}(z)=-2^{-\a-1} \{2\gamma_E+\psi([1+\a+\sigma_{\a,0}(z)]/2) 
+\psi([1+\a-\sigma_{\a,0}(z)]/2)\}, \no  \\
& \hspace*{6.8cm}  z\in\rho(T_{F,\a,0}), \; \a \in [1,\infty), \; \b=0,   \no   \\
&\sigma(T_{F,\a,0}) = \{n(n+1+\a)\}_{n\in\bbN_0},  \quad \a \in [1,\infty), \; \b = 0. 
\end{align}
\noindent 
{\boldmath $\mathbf{(III)}$ {\bf The Case $\a \in [1,\infty)$ and $\b\in(0,1)$:}}
\begin{align}
\begin{split}
\psi_{0,\a,\b}(z,x)&=2^{-\a-1}\b^{-1}y_{2,\a,\b,-1}(z,x)+y_{1,\a,\b,-1}(z,x) m_{0,\a,\b}(z),
\\
 m_{0,\a\,\b}(z)&= \b^{-1} 2^{-1 - \a - \b}\f{-\Gamma(1-\b)}{\Gamma(1+\b)}    \\
& \quad \times \f{\Gamma([1+\a+\b+\sigma_{\a,\b}(z)]/2)
\Gamma([1+\a+\b-\sigma_{\a,\b}(z)]/2)}{\Gamma([1+\a-\b+\sigma_{\a,\b}(z)]/2)
\Gamma((1+\a-\b-\sigma_{\a,\b}(z))/2)},      \\
&\hspace{3.55cm} z\in\rho(T_{F,\a,\b}), \; \a \in [1,\infty), \; \b\in(0,1),      \\
\sigma(T_{F,\a,\b}) &= \{n(n+1+\a+\b)\}_{n\in\bbN_0},   \quad \a \in [1,\infty), \; \b\in(0,1).
\end{split}
\end{align}
\noindent 
{\boldmath $\mathbf{(IV)}$ {\bf The Case $\a \in (-\infty,-1]$ and $\b\in(-1,0)$:}}  
\begin{align}
\begin{split}
\psi_{0,\a,\b}(z,x)&=y_{1,\a,\b,-1}(z,x)-2^{-\a-1}\b^{-1}y_{2,\a,\b,-1}(z,x) m_{0,\a,\b}(z), \\
 m_{0,\a\,\b}(z)&=2^{1+\a+\b}\b \f{\Gamma(1+\b)}{\Gamma(1-\b)}    \\
& \quad \times \f{\Gamma([1-\a-\b+\sigma_{\a,\b}(z)]/2)\Gamma([1-\a-\b-\sigma_{\a,\b}(z)]/2)}
{\Gamma([1+\b-\a+\sigma_{\a,\b}(z)]/2)\Gamma([1+\b-\a-\sigma_{\a,\b}(z)]/2)}, \\
&\hspace{2.55cm} \ z\in\rho(T_{F,\a,\b}), \; \a \in (-\infty, -1], \; \b\in(-1,0),  \\
\sigma(T_{F,\a,\b}) &= \{(n-\a-\b)(n+1)\}_{n\in\bbN_0},  \quad \a \in (-\infty, -1], \; \b\in(-1,0). 
\end{split}
\end{align}
\noindent 
{\boldmath $\mathbf{(V)}$ {\bf The Case $\a \in (-\infty, -1]$ and $\b=0$:}} 
\begin{align}
& \hspace{-.3cm}  \psi_{0,\a,0} (z,x)=-2^{-\a-1}y_{2,\a,0,-1}(z,x)+y_{1,\a,0,-1}(z,x)m_{0,\a,0}(z), \no \\
& m_{0,\a,0}(z)=-2^{-\a-1} \{2\gamma_E+\psi([1-\a+\sigma_{\a,0}(z)]/2)+\psi([1-\a-\sigma_{\a,0}(z)]/2)\}, 
\no \\
& \hspace*{6.25cm} z\in\rho(T_{F,\a,0}), \; \a \in (-\infty,-1], \; \b = 0,  \no \\
&\sigma(T_{F,\a,0}) = \{(n-\a)(n+1)\}_{n\in\bbN_0}, \quad \a \in (-\infty,-1], \; \b = 0.
\end{align}
\noindent 
{\boldmath $\mathbf{(VI)}$ {\bf The Case $\a \in (-\infty,-1]$ and $\b\in(0,1)$:}} 
\begin{align}
\begin{split}
\psi_{0,\a,\b} (z,x)&=2^{-\a-1}\b^{-1}y_{2,\a,\b,-1}(z,x)+y_{1,\a,\b,-1}(z,x)m_{0,\a,\b}(z), \\
m_{0,\a\,\b}(z)&= - \b^{-1} 2^{-1 - \a - \b} \f{\Gamma(1-\b)}{\Gamma(1+\b)}      \\
& \quad \times \f{\Gamma([1+\b-\a+\sigma_{\a,\b}(z)]/2)
\Gamma([1+\b-\a-\sigma_{\a,\b}(z)]/2)}{\Gamma([1-\a-\b+\sigma_{\a,\b}(z)]/2)
\Gamma([1-\a-\b-\sigma_{\a,\b}(z)]/2)},      \\
&\hspace{2.95cm}  z\in\rho(T_{F,\a,\b}), \; \a \in (-\infty,-1], \; \b\in(0,1),\\
\sigma(T_{F,\a,\b}) &= \{(n-\a)(n+1+\b)\}_{n\in\bbN_0},  \quad \a \in (-\infty, -1], \; \b\in(0,1). 
\end{split}  \lb{6.32}
\end{align}

\medskip

\appendix

\section{The Hypergeometric and Jacobi Differential Equations} \lb{sA}
\hfill

In this appendix we provide the connection between the hypergeometric differential equation (cf. \cite[Sect. 15.5]{AS72})
\begin{align}\lb{A.1}
\xi(1-\xi)\ddot{w}(\xi)+[c-(a+b+1)\xi]\dot{w}(\xi)-abw(\xi)=0, \quad \xi \in (0,1),
\end{align}
(where $\dot\ =d/d\xi$) and the Jacobi differential equation
\begin{align}\lb{A.2}
\begin{split}
\tau_{\a,\b} y(z,x) = -(1-x^2)y''(z,x)+[\a-\b+(\a+\b+2)x]y'(z,x) = z y(z,x),&   \\
\a,\b\in\bbR,\ x\in(-1,1),&
\end{split}
\end{align}
(where $' =d/dx$). Making the substitution $\xi=(1+x)/2$ in \eqref{A.2} yields
\begin{align}
\begin{split}\lb{A.3}
\xi(1-\xi)\ddot{y}(z,\xi)+[\b+1-(\a+\b+2)\xi]\dot{y}(z,\xi)+zy(z,\xi)=0,&    \\
\a,\b\in\bbR,\ \xi\in(0,1).& 
\end{split}
\end{align}
which is equal to \eqref{A.1} once one identifies, 
\begin{align}
& a = [1+ \al + \b +\sigma_{\a,\b}(z)]/2,\quad b =  [1+ \al + \b - \sigma_{\a,\b}(z)]/2,  
\quad  c= 1 + \b,    \no \\
& \sigma_{\a,\b}(z) = \big[(1+\a+\b)^2+4z\big]^{1/2}.     \lb{A.12}
\end{align}

At the endpoint $x = -1$ of the Jacobi equation the substitution used to arrive at \eqref{A.3} yields $\xi=0$, hence we next consider solutions of \eqref{A.1} near $\xi=0$ (cf. \cite[Eqs.~15.5.3, 15.5.4]{AS72}) (analogous solutions near $\xi = 1$ are found in \eqref{A.13a})
\begin{align}
& w_{1,0}(\xi) = F(a,b;c;\xi) = \sum_{n \in \bbN_0} \f{(a)_n (b)_n}{(c)_n} \f{\xi^n}{n!}, \quad 
a, b \in \bbC, \; c \in \bbC \backslash (-\bbN_0),   \no \\
& w_{2,0}(\xi) = \xi^{1-c}F(a-c+1,b-c+1;2-c;\xi),\quad a, b \in \bbC, \; (c-1) \in \bbC \backslash \bbN,   \lb{A.5} \\
& \hspace*{9.58cm} \xi \in (0,1).     \no 
\end{align}
Here $F(\dott,\dott;\dott;\dott)$ (frequently written as $\mathstrut_2F_1(\dott,\dott;\dott;\dott)$) denotes the hypergeometric function $($see, e.g., \cite[Ch.~15]{AS72}$)$, $\psi(\dott) = \Gamma'(\dott)/\Gamma(\dott)$ the Digamma function, $\gamma_{E} = - \psi(1) = 0.57721\dots$ represents Euler's constant, and 
\begin{equation}
(\zeta)_0 =1, \quad (\zeta)_n = \Gamma(\zeta + n)/\Gamma(\zeta), \; n \in \bbN, 
\quad \zeta \in \bbC \backslash (-\bbN_0), 
\end{equation}
abbreviates Pochhammer's symbol (see, e.g., \cite[Ch.~6]{AS72}). 

In addition,
\begin{equation}
\text{$w_{1,0}$ and $w_{2,0}$ are linearly independent if } \, c \in \bbC \backslash \bbZ,    
\end{equation}
which can be seen by noticing the different behaviors of $w_{1,0}(\xi)$, $w_{2,0}(\xi)$ around $\xi = 0$. 
One notes that only the case $c = 1 + \beta \in (0,2)$ is needed. Thus, for $c = 1$ we will use instead 
\begin{align}
& w_{1,0}(\xi) = F(a,b;1;\xi), \quad a, b \in \bbC,  \no \\
& w_{2,0}^{\ln}(\xi) = F(a,b;1;\xi) \, \ln(\xi)+\sum_{n\in\bbN}\f{(a)_n(b)_n}{(n!)^2}\xi^n \lb{A.8} \\
&\quad \times [\psi(a+n)-\psi(a)+\psi(b+n)-\psi(b)-2\psi(n+1)-2\gamma_E], 
\quad a,b \in \bbC \backslash (- \bbN_0),    \no \\
& \hspace*{10.67cm} \xi \in (0,1),    \no 
\end{align}
where the superscipt ``$\ln$'' indicates the presence of a logarithmic term (familiar from Frobenius theory).

Using \eqref{A.12} in formulas \eqref{A.5} and \eqref{A.8}, one obtains for the solutions of the Jacobi differential equation 
$\tau_{\a,\b} y(z,\dott) = z y(z,\dott)$ (cf.\ \eqref{A.2}) near $x=-1$, 
\begin{align}
& y_{1,\a,\b,-1}(z,x) = F(a_{\a,\b,\sigma_{\a,\b}(z)},a_{\a,\b,-\sigma_{\a,\b}(z)};1+\b;(1+x)/2),  \lb{A.13} \\
& \hspace*{7.5cm} \b \in \bbR \backslash (-\bbN),        \no \\
& y_{2,\a,\b,-1}(z,x) = (1+x)^{-\b} F(a_{\a,-\b,\sigma_{\a,\b}(z)},a_{\a,-\b,-\sigma_{\a,\b}(z)};1-\b;(1+x)/2), \no \\
& \hspace*{9.2cm} \b \in \bbR\backslash \bbN_0,   \lb{A.14} \\
& y_{2,\a,0,-1}(z,x) =F(a_{\a,0,\sigma_{\a,0}(z)},a_{\a,0,-\sigma_{\a,0}(z)};1;(1+x)/2) \, \ln((1+x)/2)     \no \\
&\quad +\sum_{n\in\bbN} \f{(a_{\a,0,\sigma_{\a,0}(z)})_n (a_{\a,0,-\sigma_{\a,0}(z)})_n}{2^n(n!)^2} 
(1+x)^n  \lb{A.15} \\
& \hspace*{1.4cm} \times [\psi(a_{\a,0,\sigma_{\a,0}(z)}+n) - \psi(a_{\a,0,\sigma_{\a,0}(z)}) 
+ \psi(a_{\a,0,-\sigma_{\a,0}(z)}+n)   \no \\
& \hspace*{1.9cm} -\psi(a_{\a,0,-\sigma_{\a,0}(z)}) -2\psi(n+1)-2\g_E], \quad \b = 0,    \no \\
& \hspace*{5.1cm} \a\in\bbR,\ z\in\bbC,\ x\in(-1,1),    \no 
\end{align}
where we abbreviated 
\begin{align}
a_{\mu,\nu,\pm\sigma} = [1 + \mu+\nu \pm \sigma]/2,\quad \mu, \nu, \sigma \in \bbC.    
\end{align}

Again one observes that for $z\in \bbC$, 
$y_{1,\a,\b,-1}(z,\dott)$ and $y_{2,\a,\b,-1}(z,\dott)$ are linearly independent for $\a \in \bbR$, $\b \in \bbR \backslash \bbZ$. 
Similarly, for $z\in \bbC$, $y_{1,\a,0,-1}(z,\dott)$ and 
$y_{2,\a,0,-1}(z,\dott)$ are linearly independent for $\a \in \bbR$.

In precisely the same manner solutions of \eqref{A.1} are given by 
\begin{align} \lb{A.13a}
    w_{1,1}(\xi) &= F(a,b;a+b-c+1;1-\xi), \quad 
a, b \in \bbC, \; c-a-b \in \bbC \backslash \bbN,    \no  
    \\  
    w_{2,1}(\xi) &= (1-\xi)^{c-a-b}F(c-a,c-b;c-a-b+1;1-\xi), 
\\\no
\;& \hspace*{3.75cm} 
a, b \in \bbC, \; a+b-c \in \bbC \backslash \bbN, 
\end{align}
and for $a+b=c$,
\begin{align}
w_{1,1}(\xi) &= F(a,b;1;1-\xi), \quad a, b \in \bbC,  \no
\\
w_{2,1}^{\ln}(\xi) &= F(a,b;1;1-\xi)\ln(1-\xi)+\sum_{n\in\bbN}\f{(a)_n(b)_n}{(n!)^2}(1-\xi)^n  \lb{B.10} \\
        &\quad\, \times[\psi(a+n)-\psi(a)+\psi(b+n)-\psi(b)-2\psi(n+1)-2\gamma_E], \no
        \\\no
\;& \hspace*{6.5cm} 
a, b \in \bbC, \; \xi \in (0,1). 
\end{align}
which are obtained from \eqref{A.5} and \eqref{A.8} by the change of variables
\begin{equation} \lb{B.1}
    (a, \, b, \, c, \, \xi) \rightarrow (a, \, b, \, a+b-c+1, \, 1 - \xi).
\end{equation}
Together with the identification $x = (1+\xi)/2$ and \eqref{A.12} one obtains the following solutions of $\tau_{\a,\b} y(z,\dott) = z y(z,\dott)$ near $x=+1$, 
\begin{align} 
& y_{1,\a,\b,+1}(z,x) = F(a_{\a,\b,\sigma_{\a,\b}(z)},a_{\a,\b,-\sigma_{\a,\b}(z)};1+\a;(1-x)/2),   \lb{A.22} \\
& \hspace{7.5cm} \a \in \bbR \backslash (- \bbN),       \no \\ 
& y_{2,\a,\b,+1}(z,x) = (1-x)^{-\a} F(a_{-\a,\b,\sigma_{\a,\b}(z)},a_{-\a,\b,-\sigma_{\a,\b}(z)};1-\a;(1-x)/2) \no \\
& \hspace*{9cm} \a \in \bbR \backslash \bbN,   \lb{A.23} \\
& y_{2,0,\b,+1}(z,x) = F(a_{0,\b,\sigma_{0,\b}(z)},a_{0,\b,-\sigma_{0,\b}(z)};1;(1-x)/2) \, \ln((1-x)/2)    \no \\
&\quad +\sum_{n\in\bbN} \f{(a_{0,\b,\sigma_{0,\b}(z)})_n (a_{0,\b,-\sigma_{0,\b}(z)})_n}{2^n(n!)^2}
(1-x)^n    \lb{A.24} \\
& \hspace*{1.4cm} \times [\psi(a_{0,\b,\sigma_{0,\b}(z)}+n) - \psi(a_{0,\b,\sigma_{0,\b}(z)})
+\psi(a_{0,\b,-\sigma_{0,\b}(z)}+n)    \no \\ 
& \hspace*{1.9cm} -\psi(a_{0,\b,-\sigma_{0,\b}(z)}) - 2\psi(n+1)-2\g_E], \quad \a =0,     \no \\
& \hspace*{5.1cm} \b\in\bbR,\ z\in\bbC,\ x\in(-1,1).    \no
\end{align}
Again, for $z\in \bbC$, $y_{1,\a,\b,+1}(z,\dott)$ and $y_{2,\a,\b,+1}(z,\dott)$ are linearly independent for $\a \in \bbR \backslash \bbZ$, $\b \in \bbR$. Similarly, for $z\in \bbC$, $y_{1,0,\b,+1}(z,\dott)$ and $y_{2,0,\b,+1}(z,\dott)$ are linearly independent for $\b \in \bbR$. 

In the limit-point case at $x = 1$, where $\a \in (-\infty,-1]\cup[1,\infty)$, one only needs the principal solutions, which are $y_{1,\a,\b,+1}(z,\dott)$ for $\a\geq 1$ and $y_{2,\a,\b,+1}(z,\dott)$ for $\a\leq -1$. Thus, one concludes from \eqref{A.22} and \eqref{A.23} that these case are already covered, and one does not have to define an additional solution for $\a \in \bbZ \setminus \lbrace 0 \rbrace$. 

Since $(a_{\a,\b,\sigma_{\a,\b}(z)})_{n} (a_{\a,\b,-\sigma_{\a,\b}(z)})_{n}$, $n \in \bbN_0$,
depends polynomially on $z \in \bbC$, one infers that 
\begin{equation}
\text{for fixed $x \in (0,1)$, $y_{j,\a,\b,\pm 1}(z,x)$, $j=1,2$, are entire with respect to $z \in \bbC$.}  
\end{equation} 
Moreover $y_{j,\a,\b,\pm 1}(z,x)$ satisfy the relations (cf.~ \eqref{A.33})
\begin{align} \lb{A.26}
y_{1,\a,\b,-1}(z,x) &= (1+x)^{-\b} y_{2,\a, -\b,-1}(z+(1+\a)\b,x),
\\&\hspace{3.15cm} \a \in \bbR, \, \b \in \bbR \setminus \lbrace 0 \rbrace,
\no\\\lb{A.27}
y_{2,\a,\b,-1}(z,x) &= (1+x)^{-\b} y_{1,\a,-\b,-1}(z+(1+\a)\b,x),
\\\no
&\hspace*{3.15cm} \a \in \bbR, \, \b \in \bbR \setminus \lbrace 0 \rbrace,
\\\no
\\\lb{A.28}
y_{1,\a,\b,+1}(z,x) &= (1-x)^{-\a} y_{2,-\a, \b,+1}(z+(1+\b)\a,x),
\\\no
&\hspace*{3.15cm}\a \in \bbR \setminus \lbrace 0 \rbrace, \, \b \in \R,
\\\lb{A.29}
y_{2,\a,\b,+1}(z,x) &= (1-x)^{-\a} y_{1,-\a,\b,+1}(z+(1+\b)\a,x),
\\\no
&\hspace*{3.17cm}\a \in \bbR \setminus \lbrace 0 \rbrace, \, \b \in \R, 
\end{align}
where we used the fact
\begin{align} \lb{A.30}
    \s_{\a,\b}(z) =
    \begin{cases}
    \s_{\a,-\b}(z+(1+\a)\b),
    \\
    \s_{-\a,\b}(z+(1+\b)\a),
    \\
    \s_{-\a,-\b}(z+\a+\b).
    \end{cases}
\end{align}
\medskip

\begin{remark}\lb{rA1}
We conclude this appendix by briefly discussing Jacobi polynomials and quasi-rational eigenfunctions. The $n$th Jacobi polynomial is defined as (see \cite[Eq. 18.5.7]{OLBC10})
\begin{align}
\begin{split}
    P_n^{\a, \b}(x) := \dfrac{(\a+1)_n}{n!} F(-n, \, n + \a + \b+1; \,
     \a + 1; \, (1-x)/2),
     \\
     n\in\N_0,\ -\a\notin\N,\ -n-\a-\b-1 \notin \N,
\end{split}
\end{align}
and can be defined by continuity for all parameters $\a, \b \in \R$. Note that $P_n^{\a,\b}(x)$ is a polynomial of degree at most $n$, and has strictly smaller degree if and only if $-n-\a-\b \in \lbrace 1, \dots , n \rbrace$ (cf.~\cite[p.~64]{Sz75}). It satisfies the equation
\begin{align}
    \tau_{\a,\b} P_n^{\a, \b}(x) = \l^{\a,\b}_n  P_n^{\a, \b}(x),
\end{align}
with 
\begin{align}
    \l^{\a,\b}_n := n(n+1+\a+\b).
\end{align}
In particular, one can verify that the Jacobi polynomials are solutions of the Jacobi differential equation \eqref{A.2} with Neumann boundary conditions at $x = +1$ (resp. $x=-1$) if $\a\in(-1,0)$ (resp. $\b\in(-1,0)$) and Friedrichs boundary conditions if $\a\geq 0$ (resp. $\b\geq0$).

More generally, all quasi-rational solutions, meaning the logarithmic derivative being rational, can be derived from the the Jacobi polynomials together with 
\begin{align}\lb{A.33}
& (1+x)^{-\b} \circ \tau_{\a,-\b} \circ (1+x)^\b = \tau_{\a,\b} + (1+\a)\b,
    \no\\
& (1-x)^{-\a} \circ \tau_{-\a,\b} \circ (1-x)^\a = \tau_{\a,\b} + (1+\b)\a,
    \\\no
& (1-x)^{-\a}(1+x)^{-\b} \circ \tau_{-\a,-\b} \circ (1-x)^\a(1+x)^\b = \tau_{\a,\b} + \a+\b,
    \end{align}
where $(1+x)^{\pm \b}$ and $(1-x)^{\pm \a}$ are regarded as formal multiplication operators. This is summarized in Table \ref{table1}, which is taken from \cite{Bo19}.
\begin{table}
    \centering
    \begin{tabular}{|l| l|}
    \hline
      Eigenfunctions   & Eigenvalues  \\
    \hline
    \hline
    $P_n^{\a,\b}(x)$ & $n(n+1+\a+\b)$
    \\
    $(1-x)^{-\a} P_n^{-\a,\b}(x)$ & $n(n+1-\a+\b)-\a(1+\b)$
    \\
    $(1+x)^{-\b} P_n^{\a,-\b}(x)$ & $n(n+1+\a-\b)-\b(1+\a)$
    \\
    $(1-x)^{-\a}(1+x)^{-\b} P_n^{-\a,-\b}(x)$ & $n(n+1-\a-\b)-(\a+\b)$
    \\
    \hline
    \end{tabular}
    \vspace{10pt}
    \caption{Formal quasi-rational eigensolutions of $\tau_{\a,\b}$}
    \label{table1}
\end{table}
Here $(1-x)^{-\a} P_n^{-\a,\b}(x)$ satisfy at $x = +1$ the Friedrichs boundary condition for $\a \leq 0$ and Neumann for $\a \in (0,1)$, while at $x = -1$ they satisfy the Friedrichs for $\b \geq 0$ and Neumann for $\b \in (-1,0)$. For $(1+x)^{-\b} P_n^{\a,-\b}(x)$ the roles of $\a$ and $\b$ interchange compared to the last case, meaning Friedrichs at $x = +1$ for $\a \geq 0$, Neumann for $\a \in (-1,0)$, and at $x = -1$, Friedrichs for $\b \leq 0$, Neumann for $\b \in (0,1)$. Finally $(1-x)^{-\a}(1+x)^{-\b} P_n^{-\a,-\b}(x)$ satisfy at $x = +1$ (resp. $x=-1$) the Friedrichs boundary condition for $\a \leq 0$ (resp. $\b\leq 0$) and Neumann for $\a \in (0,1)$ (resp. $\b\in(0,1)$).
\hfill$\diamond$
\end{remark}

\section{Connection Formulas}\lb{sB}

In this appendix we provide the connection formulas utilized to find the solution behaviors in Appendix \ref{sC}. We express them using $w_{1,0}(\xi)$ and $w_{2,0}(\xi)$ $\big(w_{2,0}^{\ln}(\xi)\big)$ and their analogs $w_{1,1}(\xi)$ and $w_{2,1}(\xi)$ $\big(w_{2,1}^{\ln}(\xi)\big)$ at the endpoint $\xi = 1$. 

We recall the relations \eqref{A.12} connecting the parameters $a,b,c$ and $\a, \b$. 
\\[1mm] 
\noindent 
$\mathbf{(I)}$ {\it The case $\a \in \R \backslash \Z$, $\b \in (-1,1) \backslash \{0\}$, that is, $c \in (0,2) \backslash \{1\}$, $a+b-c \in \R \backslash \Z$}\,: \\[1mm] 
The two connection formulas are given by (cf. \cite[Eq.~15.10.21--22]{Ov20})
\begin{align} \lb{B.4}
w_{1,0}(\xi) &= \f{\Gamma(c)\Gamma(c-a-b)}{\Gamma(c-a)\Gamma(c-b)}w_{1,1}(\xi)+\f{\Gamma(c)\Gamma(a+b-c)}{\Gamma(a)\Gamma(b)}w_{2,1}(\xi),  \\
w_{2,0}(\xi) &= \frac{\G(2-c)\G(c-a-b)}{\G(1-a)\G(1-b)}w_{1,1}(\xi)+\frac{\G(2-c)\G(a+b-c)}{\G(a-c+1)\G(b-c+1)}w_{2,1}(\xi).    \lb{B.5}
\end{align}
One notes that poles occur on the right-hand side of \eqref{B.4}, \eqref{B.5} whenever $(a+b-c) \in \bbZ$. Using \eqref{B.1} one can also express $w_{1,1}(\xi)$ or $w_{2,1}(\xi)$ as a linear combination of $w_{1,0}(\xi)$ and 
$w_{2,0}(\xi)$:
\begin{align}
    w_{1,1}(\xi) &= \f{\Gamma(a+b-c+1)\Gamma(1-c)}{\Gamma(a-c+1)\Gamma(b-c+1)}w_{1,0}(\xi)+\f{\Gamma(a+b-c+1)\Gamma(c-1)}{\Gamma(a)\Gamma(b)}w_{2,0}(\xi),     \lb{B.6} \\
    w_{2,1}(\xi) &= \frac{\G(1+c-a-b)\G(1-c)}{\G(1-a)\G(1-b)}w_{1,0}(\xi)+\frac{\G(1+c-a-b)\G(c-1)}{\G(c-a)\G(c-b)}w_{2,0}(\xi),     \lb{B.7} 
\end{align}
though, in the following we shall only write down one pair of connection formulas for brevity.
\noindent 
$\mathbf{(II)}$ {\it The case $\a = 0, \ \b \in \bbR \backslash \bbZ$, that is, $c \in \bbR \backslash \bbZ$, 
$a+b = c$}\,: \\[1mm] 
The solution  $w_{1,0}(\xi) =  F(a,b;a+b;\xi)$ can be expanded at $\xi = 1$ (cf. \cite[Eq. 15.3.10]{AS72}):
\begin{align}
\begin{split} \lb{B.8}
F(a,b;a+b;\xi)&=\f{\Gamma(a+b)}{\Gamma(a)\Gamma(b)}\sum_{n\in\N_0} \f{(a)_n(b)_n}{(n!)^2} 
[2\psi(n+1)-\psi(a+n)-\psi(b+n)\\
&\hspace*{4.05cm} -\ln(1-\xi)](1-\xi)^n.
\end{split}
\end{align}
Meanwhile, two linearly independent solutions at $\xi = 1$ are taken from \eqref{B.10}. The connection formula for $w_{1,1}(\xi)$ is given by \eqref{B.6} with $a+b=c$. To obtain a second connection formula one compares the expansion of $w_{2,1}^{\ln}(\xi)$ at $\xi = 1$ with the expansion of $F(a,b;a+b;\xi)$ at $\xi = 1$, using \eqref{B.8}, and then obtains 
\begin{align} \lb{B.11}
    w_{2,1}^{\ln}(\xi) =&  - [\psi(1-a)+\psi(1-b)+2\g_{E}] \frac{\G(1-a-b)}{\G(1-a)\G(1-b)} w_{1,0}(\xi) \no
    \\
    &- [\psi(a)+\psi(b)+2\g_{E}] \frac{\G(a+b-1)}{\G(a)\G(b)}w_{2,0}(\xi).
\end{align}
\noindent 
$\mathbf{(III)}$ {\it The case $\a \in \R \backslash \Z, \, \b = 0$, that is, $c = 1$, $a + b \in \R \backslash \Z$}\,: \\[1mm] 
This case is analogous to the previous case, with the roles of $\a$ and $\b$ interchanged. Concretely, this means that the connection formulas \eqref{B.8} and \eqref{B.11} must be changed through the renaming \eqref{B.1} with $c \rightarrow a+b-c+1=a+b$, as $c = 1$. As $c$ does not appear in \eqref{B.8} and \eqref{B.11} (it was eliminated via 
$c = a+b$), one can adopt the aforementioned formulas directly, only changing the second index in the $w$'s\footnote{Formula \eqref{B.15} could have been obtained directly from \eqref{B.4} by setting $c = 1$.}
\begin{align} 
    w_{1,0}(\xi) &= \frac{\G(1-a-b)}{\G(1-a)\G(1-b)}w_{1,1}(\xi)+ \frac{\G(a+b-1)}{\G(a)\G(b)}w_{2,1}(\xi),
   \lb{B.15} \\
    w_{2,0}^{\ln}(\xi) &=- [\psi(1-a)+\psi(1-b)+2\g_{E}] \frac{\G(1-a-b)}{\G(1-a)\G(1-b)}w_{1,1}(\xi)
   \no \\
    & \quad \, - [\psi(a)+\psi(b)+2\g_E] \frac{\G(a+b-1)}{\G(a)\G(b)}w_{2,1}(\xi).
    \lb{B.16} 
\end{align}    
\noindent 
$\mathbf{(IV)}$ {\it The case $\a = \b = 0$, that is, $ a + b = c = 1$}\,: \\[1mm] 
For $\a = 0$ and $\b = 0$ the Jacobi differential expression \eqref{6.1} becomes the Legendre differential expression. Since this case was treated in detail in \cite{GLN20}, we shall only present the connection formulas for completeness.  

The special solutions $w_{1,i}(\xi)$ and $w_{2,i}^{\ln}(\xi)$ for $i = 1,2$ are given by \eqref{A.8} and \eqref{B.10}, respectively. Note that the following relations hold
\begin{align} \lb{B.19}
    w_{1,1}(\xi) = w_{1,0}(1-\xi), \qquad
    w_{2,1}^{\ln}(\xi) = w_{2,0}^{\ln}(1-\xi).
\end{align}
Using \cite[Eq.~15.3.10]{AS72} together with $w_{1,0} = F(a,b;a+b, \xi)$ 
and Euler's famous reflection formula, $\G(z)\G(1-z) = \pi \csc(\pi z)$ (cf. \cite[Eq.~6.1.17]{AS72}),
one obtains 
\begin{align}
    w_{1,0}(\xi) = - \pi^{-1} \sin(\pi a) \big([\psi(a)+\psi(b)+2\g_E]w_{1,1}(\xi)+w_{2,1}^{\ln}(\xi)\big).
\end{align}
The two relations \eqref{B.19} immediately imply 
\begin{align}\label{B.22}
    w_{1,1}(\xi) &= - \pi^{-1} \sin(\pi a) \big([\psi(a)+\psi(b)+2\g_E]w_{1,0}(\xi)+w_{2,0}^{\ln}(\xi)\big),
    \\
\label{B.24}
    w_{2,1}^{\ln}(\xi) &=  \pi^{-1} \sin(\pi a) \big[\big([\psi(a)+\psi(b)+2\g_E]^2-\pi^2 [\sin(\pi a)]^{-2}\big) 
     w_{1,0}(\xi)
    \no\\
    &\quad + [\psi(a)+\psi(b)+2\g_E] w_{2,0}^{\ln}(\xi)\big].
\end{align}

\section{Behavior of $y_{j,\a,\b,\mp 1}(z,x)$, $j=1,2$, near $x=\pm 1$} \lb{sC}

In this appendix we focus on the generalized boundary values for the solutions $y_{j,\a,\b,-1}(z,x),\ j=1,2$ at $x= \mp 1$. One obtains for $z\in\bbC$,
\begin{align}\lb{C.13}
\begin{split}
\wti y_{1,\a,\b,-1}(z,-1)&=\begin{cases}
1, & \b\in(-1,0),\\
0, & \b=0,\\
0, & \b\in(0,1),
\end{cases}\\
\wti y_{1,\a,\b,-1}^{\, \prime}(z,-1)&=\begin{cases}
0, & \b\in(-1,0),\\
1, & \b=0,\\
1, & \b\in(0,1),
\end{cases}\\
\wti y_{2,\a,\b,-1}(z,-1)&=\begin{cases}
0, & \b\in(-1,0),\\
-2^{\a+1}, & \b=0,\\
\b 2^{\a+1}, & \b\in(0,1),
\end{cases}\\
\wti y_{2,\a,\b,-1}^{\, \prime}(z,-1)&=\begin{cases}
-\b 2^{\a+1}, & \b\in(-1,0),\\
0, & \b=0,\\
0, & \b\in(0,1),
\end{cases}
\end{split}
\quad \a\in\bbR,
\end{align}
and employing connection formulas for the endpoint $x=+1$,
\begin{align}
\no \begin{split}
\wti y_{1,\a,\b,-1}(z,1)&=\begin{cases}
\dfrac{\Gamma(1+\b)\Gamma(-\a)}{\Gamma(a_{-\a,\b,\sigma_{\a,\b}(z)})\Gamma(a_{-\a,\b,-\sigma_{\a,\b}(z)})}, & \a\in(-1,0),\\[3mm]
\dfrac{- 2^{1+\a+\b}\Gamma(1+\a)\Gamma(1+\b)}{\Gamma(a_{\a,\b,\sigma_{\a,\b}(z)})\Gamma(a_{\a,\b,-\sigma_{\a,\b}(z)})}, & \a\in[0,1),
\end{cases}
\end{split}\\[1mm]
\no \begin{split}
\wti y_{1,\a,\b,-1}^{\, \prime}(z,1)&=\begin{cases}
\dfrac{2^{1+\a+\b}\Gamma(1+\a)\Gamma(1+\b)}{\Gamma(a_{\a,\b,\sigma_{\a,\b}(z)})\Gamma(a_{\a,\b,-\sigma_{\a,\b}(z)})}, & \a\in(-1,0),\\[3mm]
\dfrac{-\Gamma(1+\b)}{\Gamma(a_{0,\b,\sigma_{0,\b}(z)})\Gamma(a_{0,\b,-\sigma_{0,\b}(z)})}[2\gamma_E\\
\quad+\psi(a_{0,\b,\sigma_{0,\b}(z)})+\psi(a_{0,\b,-\sigma_{0,\b}(z)})], & \a=0,\\[3mm]
\dfrac{\Gamma(1+\b)\Gamma(-\a)}{\Gamma(a_{-\a,\b,\sigma_{\a,\b}(z)})\Gamma(a_{-\a,\b,-\sigma_{\a,\b}(z)})}, & \a\in(0,1),
\end{cases}
\end{split}\\[1mm]
&\lb{C.14} \hspace{6.5cm} \b\in(-1,1),\\[1mm]
\no\begin{split}
\wti y_{2,\a,\b,-1}(z,1)&=\begin{cases}
\dfrac{2^{-\b}\Gamma(1-\b)\Gamma(-\a)}{\Gamma(a_{-\a,-\b,\sigma_{\a,\b}(z)})\Gamma(a_{-\a,-\b,-\sigma_{\a,\b}(z)})}, & \a\in(-1,0),\\[3mm]
\dfrac{-2^{\a+1}\Gamma(1+\a)\Gamma(1-\b)}{\Gamma(a_{\a,-\b,\sigma_{\a,\b}(z)})\Gamma(a_{\a,-\b,-\sigma_{\a,\b}(z)})}, & \a\in[0,1),
\end{cases}
\end{split}\\[1mm]
\no \begin{split}
\wti y_{2,\a,\b,-1}^{\, \prime}(z,1)&=\begin{cases}
\dfrac{2^{\a+1}\Gamma(1+\a)\Gamma(1-\b)}{\Gamma(a_{\a,-\b,\sigma_{\a,\b}(z)})\Gamma(a_{\a,-\b,-\sigma_{\a,\b}(z)})}, & \a\in(-1,0),\\[3mm]
\dfrac{-2^{-\b}\Gamma(1-\b)}{\Gamma(a_{0,-\b,\sigma_{0,\b}(z)})\Gamma(a_{0,-\b,-\sigma_{0,\b}(z)})}[2\gamma_E\\
\quad+\psi(a_{0,-\b,\sigma_{0,\b}(z)})+\psi(a_{0,-\b,-\sigma_{0,\b}(z)})], & \a=0,\\[3mm]
\dfrac{2^{-\b}\Gamma(1-\b)\Gamma(-\a)}{\Gamma(a_{-\a,-\b,\sigma_{\a,\b}(z)})\Gamma(a_{-\a,-\b,-\sigma_{\a,\b}(z)})}, & \a\in(0,1),
\end{cases}
\end{split}  \\[1mm]
& \hspace{6.25cm} \b\in(-1,1)\backslash\{0\},   \lb{C.15} \\[1mm]
\no \begin{split}
\wti y_{2,\a,0,-1}(z,1)&=\begin{cases}
\dfrac{-[2\gamma_E+\psi(a_{-\a,0,\sigma_{\a,0}(z)})+\psi(a_{-\a,0,-\sigma_{\a,0}(z)})] 
\Gamma(-\a)}{\Gamma(a_{-\a,0,\sigma_{\a,0}(z)})\Gamma(a_{-\a,0,-\sigma_{\a,0}(z)})},\\[1mm]
\hspace{6.2cm}\a\in(-1,0),\\[3mm]
\dfrac{[2\gamma_E+\psi(a_{\a,0,\sigma_{\a,0}(z)})+\psi(a_{\a,0,-\sigma_{\a,0}(z)})]\Gamma(1+\a)}{2^{-\a-1}\Gamma(a_{\a,0,\sigma_{\a,0}(z)})\Gamma(a_{\a,0,-\sigma_{\a,0}(z)})},\\[1mm]
\hspace{6.1cm}\a\in[0,1),
\end{cases}
\end{split}\\[1mm]
\begin{split} 
\wti y_{2,\a,0,-1}^{\, \prime}(z,1)&=\begin{cases}
\dfrac{[2\gamma_E+\psi(a_{\a,0,\sigma_{\a,0}(z)})+\psi(a_{\a,0,-\sigma_{\a,0}(z)})]\Gamma(1+\a)}{-2^{-\a-1}\Gamma(a_{\a,0,\sigma_{\a,0}(z)})\Gamma(a_{\a,0,-\sigma_{\a,0}(z)})},\\[1mm]
\hspace{5.8cm}\a\in(-1,0),\\[1mm]
-\Gamma((1+\sigma_{0,0}(z))/2)\Gamma((1-\sigma_{0,0}(z))/2)\\
+\dfrac{[2\gamma_E
+\psi((1+\sigma_{0,0}(z))/2)+\psi((1-\sigma_{0,0}(z))/2)]^2}{\Gamma((1+\sigma_{0,0}(z))/2)
\Gamma((1-\sigma_{0,0}(z))/2)},\quad \a=0,\\[3mm]
\dfrac{-[2\gamma_E+\psi(a_{-\a,0,\sigma_{\a,0}(z)})+\psi(a_{-\a,0,-\sigma_{\a,0}(z)})]\Gamma(-\a)}
{\Gamma(a_{-\a,0,\sigma_{\a,0}(z)})\Gamma(a_{-\a,0,-\sigma_{\a,0}(z)})},\\[1mm]
\hspace{6.45cm}\a\in(0,1),  
\end{cases}   
\end{split} \no \\[1mm] 
&\hspace{7.85cm} \b = 0.     \lb{C.16}   
\end{align}


 
\end{document}